\documentclass[11pt]{article}

\usepackage[a4paper,margin=1in]{geometry}
\usepackage{amsmath,amssymb,amsthm,mathtools}
\usepackage{enumitem}
\usepackage{hyperref}
\usepackage{xcolor}
\usepackage{needspace}
\numberwithin{equation}{section}
\allowdisplaybreaks
\newtheorem{theorem}{Theorem}[section]

\newtheorem{lemma}[theorem]{Lemma}
\newtheorem{remark}[theorem]{Remark}

\title{
Decay Rate toward the Rarefaction Wave and the Viscous Contact Wave for Jin--Xin Relaxation Model
}

\author{Mutong He\thanks{Academy of Mathematics and Systems Science, Chinese Academy of Sciences, Beijing 100190, PR China. E-mail: hemutong25@mails.ucas.ac.cn},
Feimin Huang\thanks{Academy of Mathematics and Systems Science, Chinese Academy of Sciences, Beijing 100190, PR China. E-mail: fhuang@amt.ac.cn},
Tianyi Wang\thanks{School of Mathematics and Statistics, Huazhong University of Science and Technology, Wuhan 430074, PR China. E-mail: {tian-yi.wang@hust.edu.cn}},
Lingda Xu\thanks{Hetao Institute of Mathematics and Interdisciplinary Sciences, Shenzhen 518048, PR China. E-mail: xuldmath@gmail.com}
}
\date{}
\begin{document}
\maketitle

\begin{abstract}
\begingroup
This paper establishes the time-decay rates toward the rarefaction wave and the viscous contact wave for the one-dimensional Jin--Xin relaxation system under the $L^1$ integrability condition of the initial perturbation. In the case of a strictly convex flux, we employ a decomposition method to derive a logarithmic $L^1$ estimate related to a conjecture proposed by Liu\cite{liu2001lp}. Furthermore, a time-weighted $L^2$ method yields the optimal decay rate toward the rarefaction wave. Regarding the case where the flux is partially linearly degenerate, by employing a similar \(L^1\) estimate and time-weighted \(L^2\) method, along with a refined estimate of the nonlinear wave interactions, we obtain the corresponding decay rate toward the composite wave. The decay estimates are obtained without smallness assumptions on the initial perturbations or the wave strength.
\endgroup
\end{abstract}



\tableofcontents

\section{Introduction and Main Results}

This paper is concerned with the one-dimensional Jin--Xin relaxation system
\begin{equation}
\label{Original equation: Jin-Xin system}
\begin{cases}
        u_t+v_x=0,\\
        v_t+u_x=f(u)-v,
\end{cases}
        \qquad
        (t,x)\in \mathbb R_+\times \mathbb R,
\end{equation}
\noindent with initial data
\begin{equation}
\label{initial values}
    (u,v)(0,x)=(u_0(x),v_0(x)),
\end{equation}

\noindent which satisfies
\begin{equation}
        (u,v)(t,x)\to (u_\pm,f(u_\pm)),
        \quad\mathrm{as} \quad x\to\pm\infty,
\label{eq:farfield}
\end{equation}
\noindent where $u_\pm$ are two given far-field states satisfying $u_-<u_+$. Furthermore, we assume that the nonlinear flux function $f(u)$ is a given $C^1(\mathbb{R})$ function and meet the sub-characteristic condition:
\begin{equation}
    \label{sub-characteristic}
    k:=\sup\limits _{u \in M}\left|f^{\prime}(u)\right|<1, 
\end{equation}
\noindent where the state space \(M\) will be specified later.

We are interested in the asymptotic behavior  of the global-in-time solution of the Cauchy problem (\ref{Original equation: Jin-Xin system})-(\ref{initial values}), especially the influence of the properties of the nonlinear function \(f\left( u\right)\) on the asymptotic behavior. In particular, we consider that the flux \(f\) is strictly convex, that is 
\begin{equation}
\label{strictly convex}
    f''(u)\geq \kappa_0>0, \quad u\in \mathbb{R},
\end{equation}
\noindent or  the flux is smooth and genuinely nonlinear on the \(\mathbb{R}\) except for a finite interval
$I:=(a,b)\subset \mathbb{R}$, which is linearly degenerate on \(I\), that is
\begin{equation}
\label{Partially linearly degenerate condition}
    \begin{cases}f^{\prime \prime}(u)>0, &u \in(-\infty, a] \cup[b,+\infty), \\ f^{\prime \prime}(u)=0,  &u \in(a, b),\end{cases}
\end{equation}
\noindent Here the second derivatives at $a$ and $b$ are understood from the nonlinear sides, and the propagation speed is ${f}^{\prime }(u)=\bar{\lambda }$ on $[a,b]$.

Note that the relaxation system (\ref{Original equation: Jin-Xin system}) is equivalent to the following perturbed viscous conservation laws:
\begin{equation}
\label{perturbed viscous conservation laws}
    u_t+f(u)_x=u_{x x}-u_{t t},
\end{equation}
\noindent which can be easily derived by eliminating $v$ from (\ref{Original equation: Jin-Xin system}). We treat the term $u_{t t}$ as a higher-order perturbation term in (\ref{perturbed viscous conservation laws}) and expect that the long time behavior of the solutions to (\ref{Original equation: Jin-Xin system}) is similar to that for the viscous conservation laws
\begin{equation}
u_t+f(u)_x=u_{x x} .
\end{equation}
\noindent For the Riemann problem for the conservation law,
\begin{equation*}
\label{inviscid Burgers}
    \left\{\begin{array}{l}
u_t+f(u)_x=0 ,\\
\left.u\right|_{t=0}=u_0^R(x)=\left\{\begin{array}{ll}
u_{-},  \quad &x<0,  \\
u_{+},  \quad &x>0.
\end{array}\right.
\end{array}\right.
\end{equation*}

 \noindent If the flux is strictly convex, 
the entropy solution $u^r\left(\frac{x}{t}\right)$ is called the rarefaction wave, which is given explicitly as
\begin{equation}
    u^r\left(\frac{x}{t}\right)= \begin{cases}u_{-},  & x \leq f^{\prime}\left(u_{-}\right) t, \\ \left(f^{\prime}\right)^{-1}(x / t),  & f^{\prime}\left(u_{-}\right) t<x<f^{\prime}\left(u_{+}\right) t, \\ u_{+},  & x \geq f^{\prime}\left(u_{+}\right) t,\end{cases}
\end{equation}

For the asymptotic stability of the basic wave patterns (i.e., the rarefaction wave, shock wave, and contact discontinuity) for the viscous conservation law, many results have been established, see,  \cite{il1964behavior}, \cite{nishihara2001asymptotic}, \cite{liu1998behaviors}, \cite{matsumura2012asymptotic}, \cite{huang2025time}, \cite{ii1960asymptotic}, \cite{kawashima1985asymptotic}, \cite{nishihara1985note}, \cite{matsumura1994asymptotic}, \cite{kawashima1994stability}, \cite{mei1995stability}. It is also interesting to study the corresponding time-decay rate. Regarding the rarefaction wave, under the \(L^1\) integrability condition on the initial data, the time-decay rate can be obtained via the \(L^p\) energy method and the interpolation inequality, see \cite{Escobedo1991}, \cite{Kawashima2004}, \cite{Nishikawa2000}, \cite{Huang2021}. Recently, Yoshida \cite{yoshida2014decay} obtained the decay rate toward the composite wave consisting of the rarefaction wave and the viscous contact wave by applying a similar method that removes the \(L^1\) integrability assumption. As a further extension to the multi-dimensional case, Huang-Xu \cite{HuangXu2024} showed that the planar rarefaction wave is nonlinearly stable with a time decay rate without the \(L^1\) integrability condition. As for the viscous shock wave, Kawashima and Matsumura \cite{kawashima1985asymptotic} first obtained the decay rate for the scalar viscous conservation law by a weighted energy method. There have been several works on the decay properties toward viscous shocks, cf. \cite{Jones1993}, \cite{matsumura1994asymptotic}, \cite{nishihara1985note}, in which all the decay rates are obtained by time-weighted estimates. In recent years, Huang and Xu \cite{huang2022decay} obtained the time-decay rate via an \(L^p\) estimate and an area inequality (cf. \cite{dong2021large}) under the assumption of small initial perturbations without any additional \(L^1\) integrability condition.

Concerning the asymptotic behaviour of the solution to the Jin-Xin relaxation model \eqref{Original equation: Jin-Xin system}, for the case that the flux is strictly convex, the nonlinear stability of the rarefaction wave for the 1-D Jin-Xin relaxation system was originally considered by Liu \cite{liu1987hyperbolic}. Later, Zhao \cite{zhao2000nonlinear} obtained the same stability without the smallness conditions on the wave strength and the initial perturbation. Furthermore, under the condition that the second-order or third-order derivatives of the fluxes are uniformly bounded, Zou, Zhao, and Wang \cite{2012The} showed the stability under large initial perturbations, based on the elementary energy method. Corresponding results for the stability of the planar rarefaction wave in two dimensions can be found in \cite{luo1997asymptotic}, \cite{zhao2000nonlinear}, and \cite{2012The}. However, there are few results for the decay rate in time toward the rarefaction wave. In 2001, Liu \cite{liu2001lp} tried to get the decay rate by employing the \(L^1\)-contraction and the time-weighted \(L^2\) approach with \(L^1\) integrability. Nevertheless, the standard \(L^1\) stability estimate for Riemann invariants, see, e.g. \cite{natalini1996convergence}, cannot lead to the correct \(L^1\) estimate he conjectured. Liu obtained decay estimates toward relaxation rarefaction profiles; the present argument concerns the Euler rarefaction wave. As for the case that the flux is partially linearly degenerate, for the equivalent damped wave equation with a nonlinear convection term derived from the 1-D Jin-Xin relaxation system:

\begin{equation}
\label{initial wave equ}
    u_{t t}+u_t+f(u)_x=\mu u_{x x}, \quad t>0, \quad x \in \mathbb{R},  \\
\end{equation}
\noindent with initial data
\begin{equation}
\label{initial value of wave equ}
\left\{\begin{array}{l}
u(0, x)=u_0(x) \rightarrow u_{ \pm},\quad \mbox{as}\quad x \rightarrow \pm \infty, \\
u_t(0, x)=u_1(x), 
\end{array}\right.
\end{equation}

\noindent Yoshida \cite{MR4250099} first studied the stability of the composite wave consisting of the rarefaction wave and the viscous contact wave via a technical energy estimate under the special assumption \(\bar{\lambda}:=f'(0)=0\) and the smallness condition on the wave strength and the initial perturbation, where the viscous contact wave connecting \(u_{-}\) and \(u_{+}\ (u_{\pm} \in [a,b])\) is given by the exact solution of the linear convection heat equation
\begin{equation}
    \label{vc equ}
u_t+\bar{\lambda} u_x=\mu  u_{xx},\qquad t>0,\quad x \in \mathbb{R},
\end{equation}
\noindent its form is
\begin{equation}
\label{Expression of VC}
   U=U(t,x)= U\left(\frac{x-\bar\lambda t}{\sqrt{t}} ; u_{-}, u_{+}\right):=u_{-}+\frac{u_{+}-u_{-}}{\sqrt{\pi}} \int_{-\infty}^{\frac{x-\bar\lambda t}{\sqrt{4 \mu t}}} \mathrm{e}^{-\xi^2} \mathrm{~d} \xi ,
\end{equation}
\noindent where \(\mu\) represents the viscous coefficient of the viscous contact wave. Later, He, Huang, and Wang \cite{MR4983541} obtained the general stability results in which the assumptions \(\bar{\lambda}=0\) and the smallness conditions were removed. It is worth noting that the authors used a refined construction of the viscous contact wave accounting for a non-normalizable propagation speed; that is, the authors set the coefficient \(\mu\) in \eqref{vc equ} as \(1-\bar{\lambda}^2\), whose positivity is guaranteed by the sub-characteristic condition \eqref{sub-characteristic}. For the one-dimensional case, Huang, Pan, and Wang \cite{Huang2008Stability} constructed a contact wave for 1-D Jin-Xin relaxation systems, which is shown to be nonlinearly stable under small initial perturbations. The time-decay rate was also obtained by weighted energy estimates. Moreover, for the half-space problem of \eqref{initial wave equ}, Fan, Liu, and Zhao \cite{fan2013one} showed that the global solution of the half-space problem converges to the boundary layer and the decay rate was obtained by the time-space weighted estimate method, in which the initial data belongs to a weighted Sobolev space.

In this paper, we investigate the asymptotic stability and the decay rates of the perturbation for the Jin–Xin relaxation system. By employing a refined decomposition method, we obtain the \(L^1\) bound for the nonlinear flux, thereby verifying Liu's conjecture. Furthermore, by applying the time-weighted energy approach, we establish the decay rates in the following two distinct regimes. In the case of a strictly convex flux, we derive the decay rate \((1+t)^{-1/2+\delta}\) toward the rarefaction wave. Notably, by comparing this with the decay rate of the corresponding damped wave equation, we find that this is the optimal decay rate. For the partially linearly degenerate flux, we derive the decay rate \((1+t)^{-1/3+\delta}\) toward the composite wave of a rarefaction wave and a viscous contact wave, where \(\delta>0\) is an arbitrary constant. In addition, we provide the interaction estimates of the nonlinear waves for the latter case.

The key points of our work are as follows: First, for the case that the flux is strictly convex, in contrast to Liu \cite{liu2001lp}, who utilized the standard \(L^1\) stability estimate for Riemann invariants, we employ an exquisite decomposition method, inspired by a careful observation of the structure of the one-dimensional Jin--Xin relaxation system, to establish the \(L^1\) estimate, thereby obtaining logarithmic control of the perturbation. Second, under the \(L^1\) estimate we have obtained, we derive the decay rate in time toward the rarefaction wave by employing the time-weighted \(L^2\) estimate, rather than the time-space weighted estimate as in \cite{yin2008convergence}, \cite{fan2013one}, and related works. Third, regarding the case that the flux is partially linearly degenerate, by applying a similar decomposition method, we also deduce the \(L^1\) bound for the perturbation, in which we use the precise estimate of the interaction of the nonlinear waves. Finally, since \(f''\) is discontinuous in the partially linearly degenerate case, the classical integration by parts is no longer applicable when deriving higher-order derivative decay estimates. To obtain a sharper decay rate, by dividing the domain into distinct regions based on both the zeros of the smooth composite wave and the discontinuities of \(f''\), and thereby making a more refined estimate, we derive the decay rate toward the composite wave consisting of the rarefaction wave and the viscous contact wave.

\section{Preliminaries}
In this section, we provide some preliminaries necessary for the proofs of the main theorems. We now state the interpolation inequalities, which will be frequently applied later.
\begin{lemma}[Interpolation inequality]
\label{Interpolation inequality}

(1) (Nash inequality) 
If \(u\in H^1(\mathbb{R})\cap L^1(\mathbb{R})\), then
\begin{equation*}
    \label{Nash inequ}
    \|u\|_{L^2} \leq C\|u_x\|_{L^2{}}^{\frac{1}{3}}\|u\|_{L^1{}}^{\frac{2}{3}}.
\end{equation*} 

(2) If \(u\in H^2(\mathbb{R})\), then
\begin{equation}
    \label{derivative inequ}
    \|u_x\|_{L^2{}} \leq C\|u_{xx}\|_{L^2{}}^{\frac{1}{2}}\|u\|_{L^2{}}^{\frac{1}{2}}.
\end{equation}
\end{lemma}

\noindent Next, we introduce the global existence and uniform boundedness results established in \cite{zhao2000nonlinear,natalini1996convergence}.
\begin{lemma}[Global existence and uniform bounds \cite{zhao2000nonlinear,natalini1996convergence}]
\label{Global existence}
Suppose that $f\in C^1(\mathbb{R})\subset \mathrm{Lip}_{loc}(\mathbb{R})$ and $(u_0,v_0)\in L^\infty(\mathbb{R})$, and set
\[
N_0=\max\{\|u_0\|_{L^\infty},\|v_0\|_{L^\infty}\},\qquad
F(N_0)=\sup_{|u|\leq N_0}|f(u)|,\qquad B(N_0)=2N_0+F(2N_0).
\]
Then the Cauchy problem (\ref{Original equation: Jin-Xin system})-(\ref{initial values}) admits a unique globally bounded weak solution $(u,v)$, which satisfies, for $t\geq0$,
\[
\|u(t)\|_{L^\infty}\leq2N_0,\qquad\|v(t)\|_{L^\infty}\leq B(N_0).
\]
If, in addition, the initial data satisfy the regularity assumptions in \cite{zhao2000nonlinear}, then this solution is the unique global smooth solution considered below. As a direct corollary, the solution satisfies the sub-characteristic condition (\ref{sub-characteristic}) for each $t\geq0$ and $x\in\mathbb R$.
\end{lemma}
\begin{proof}[Proof of the uniform bounds]
Put $A=2N_0=2\max \left\{\left\|u_0(x)\right\|_{L^{\infty}},\left\|v_0(x)\right\|_{L^{\infty}}\right\}$, $p=u+v$ and $q=u-v$. Then
\begin{equation*}
\begin{cases}
p_t+p_x=f\bigl((p+q)/2\bigr)-(p-q)/2,\\
q_t-q_x=-f\bigl((p+q)/2\bigr)+(p-q)/2.
\end{cases}
\end{equation*}
For the comparison argument, extend $f$ outside $M=[-A,A]$ by its tangent lines at the endpoints, keeping $|f'|\leq k<1$. The off-diagonal derivatives of the two source terms are $(1+f')/2$ and $(1-f')/2$, respectively, so this system is cooperative. The spatially constant solutions
\begin{equation*}
p^\pm(t)=\pm A+f(\pm A)(1-e^{-t}),\qquad
q^\pm(t)=\pm A-f(\pm A)(1-e^{-t})
\end{equation*}
\noindent satisfy $(p^-(0),q^-(0))\leq(p(0,x),q(0,x))\leq(p^+(0),q^+(0))$ componentwise. The comparison principle therefore gives $p^-\leq p\leq p^+$ and $q^-\leq q\leq q^+$. Consequently, using the definition of $(p,q)$, one could have
\begin{equation*}
\begin{aligned}
& |u(t, x)| \leq A=2 N_0, \\
& \left|v(t, x)-\frac{f(A)+f(-A)}{2}\left(1-e^{-t}\right)\right| \leq A,
\end{aligned}
\end{equation*}
\noindent which leads to
\begin{equation*}
    |v(t, x)| \leq A+\frac{|f(A)+f(-A)|}{2} \leq A+F(A)=B\left(N_0\right),
\end{equation*}
\noindent where $F\left(N_0\right)=\sup\limits _{|u| \leq N_0}\left|f(u)\right|$. Since $u$ stays in $M$, the extension of $f$ does not change the solution of the original system. This also explains why the state interval in the sub-characteristic condition is $[-2N_0,2N_0]$, whereas the bound for $v$ involves the size of $f$ itself.
\end{proof}

\section{The strictly convex case}
This section is devoted to the proof of Theorem \ref{Main Th 1}. We consider the case where the flux function \(f\) is strictly convex, that is, \(f\) satisfies (\ref{strictly convex}).

\subsection{The smooth rarefaction wave of Type 1}
For the strictly convex case, we first introduce the construction of the smooth rarefaction wave of Type 1. Considering the Riemann problem of the inviscid Burgers equation
\begingroup
\predisplaypenalty=0
\begin{align}
    \label{1}
    \left\{
    \begin{array}{l}
        w_t+w w_x=0,\\
        w(0,x)=
        \begin{cases}
            w_-, & x<0,\\
            w_+, & x>0.
        \end{cases}
    \end{array}
    \right.
\end{align}
\endgroup
\noindent When \({w}_{ - } < {w}_{ + }\), the self-similar solution \({w}^{r}\left( {\frac{x}{t};{w}_{ - },{w}_{ + }}\right)\) of (\ref{1}) is called the centered rarefaction wave
\begin{equation*}
   w^r\left(\frac{x}{t};w_{-}, w_{+}\right)= \begin{cases}w_{-},  & x \leq w_{-} t,  \\ \frac{x}{t},  & w_{-} t<x<w_{+} t , \\ w_{+},  & x \geq w_{+} t.\end{cases}
\end{equation*}
\noindent If the flux function \(f\) is strictly convex and \({u}_{ - } < {u}_{ + }\), then the Riemann problem of the general conservation law 
\begin{equation}
    \left\{\begin{array}{l}
u_t+f(u)_x=0,  \\
\left.u\right|_{t=0}=u_0^R(x)=\left\{\begin{array}{ll}
u_{-}, \quad & x<0 , \\
u_{+},  \quad & x>0.
\end{array}\right.
\end{array}\right.
\end{equation}
\noindent has the following rarefaction wave solution
\begin{equation}
\label{rare}
u^r\left(\frac{x}{t}, u_{-},  u_{+}\right)=(\lambda)^{-1}\left(w^r\left(\frac{x}{t} ; \lambda_{-},  \lambda_{+}\right)\right)
  = \begin{cases}u_{-},  & x \leq f^{\prime}\left(u_{-}\right) t , \\ \left(\lambda\right)^{-1}(x / t),  & f^{\prime}\left(u_{-}\right) t<x<f^{\prime}\left(u_{+}\right) t,  \\ u_{+},  & x \geq f^{\prime}\left(u_{+}\right) t,\end{cases}
\end{equation}
\noindent where the characteristic speed \(\lambda \left( u\right)  = {f}^{\prime }\left( u\right)\).

Due to the singularity of \({u}^{r}\left( {\frac{x}{t};{u}_{ - },{u}_{ + }}\right)\) at \(t = 0\), which cannot serve directly as the background state, we construct the smooth rarefaction wave solution \(U^r\left( {t,x}\right)\) to approximate the rarefaction wave \({u}^{r}\left( \frac{x}{t}\right)\) of Riemann problem.

For $w^r(\frac{x}{t};w_{-}, w_{+})$, we construct corresponding smooth rarefaction wave $W(t, x)$ as
\begin{equation}
\label{2}
    \left\{\begin{array}{l}
W_t+W W_x=0,  \\
W(x,  0)=W_0(x)=\frac{1}{2}\left(w_{+}+w_{-}\right)+\widetilde{w} K_q \int_0^{\epsilon x}\left(1+y^2\right)^{-q} d y.
\end{array}\right.
\end{equation}
\noindent where $\widetilde{w}=\frac{w_{+}-w_{-}}{2}>0$ and $q>\frac{1}{2}$, $ \varepsilon>0$ are constants, $ K_q$ is the fixed constant satisfying
\begin{equation*}
K_q=\frac{1}{ \int_{0}^{+\infty}\left(1+y^2\right)^{-q} d y}.    
\end{equation*}
And, the smooth approximation of $u^r(\frac{x}{t})$ is  $U^r(t, x)=\left(f^{\prime}\right)^{-1}(W(t, x))$. Denote $\widetilde{u}=\frac{u_{+}-u_{-}}{2}>0$ ($u_+>u_-$) and let $q=1$, $\varepsilon=\widetilde{w}$, then we could have following properties of $U^r(t, x)$.

\begin{lemma}[Properties of the smooth rarefaction wave of Type 1]
\label{Properties of the smooth rarefaction wave 1}
    Let $f(u)\in C^4( \mathbb{R})$,   $f^{\prime \prime}(u)>0$, then
    
(1) $U^r(t, x)$ is the $C^3$ global solution in time of the following Cauchy problem
\begin{equation*}
    \left\{\begin{array}{l}
U^r_t+f(U^r)_x=0,  \\
\left.U^r\right|_{t=0}=U^r_0(x)=\left(f^{\prime}\right)^{-1}\left(W_0(x)\right) .
\end{array}\right.
\end{equation*}

(2) $U^r_x>0$, $u_-<U^r(t, x)<u_+$.

(3)For $ 1\leq p\leq +\infty$,  there exists a constant $C_p>0$ depending on $p$ such that for $t\geq 0$,
$$
\begin{aligned}
& \left\|U_x^r(t)\right\|_{L^p},\left\|U_t^r(t)\right\|_{L^p} \leq C_p {} h(\widetilde{u})(1+t)^{-1+\frac{1}{p}}, \\
& \left\|U_{x x}^r(t)\right\|_{L^p},\left\|U_{t t}^r(t)\right\|_{L^p},\left\|U_{t x}^r(t)\right\|_{L^p} \leq C_p {} h(\widetilde{u})(1+t)^{-\frac{3}{2}+\frac{1}{2 p}}, \\
& \left\|U_{x x x}^r(t)\right\|_{L^p},\left\|U_{t t x}^r(t)\right\|_{L^p} \leq C_p {} h(\widetilde{u})(1+t)^{-2+\frac{1}{p}}.
\end{aligned}
$$

\noindent where $h(\widetilde{u})$ is a nondecreasing function of $\widetilde{u}$ with $ \lim\limits _{\widetilde{u} \to 0} h(\widetilde{u})=0$.

(4) $\sup\limits _{x \in \mathbb{R}}\left|U^r(t,  x)-u^r\left(\frac{x}{t}\right)\right| \rightarrow 0$ as $t \rightarrow +\infty $.

(5){For $1< p\leq+\infty$, there is a constant $C_{p,\delta}>0$ such that for $t>0$
\begin{equation*}
   \|U^r(t,x)-u^r(\frac{x}{t})\|_{L^p} \leq C_{p, \delta} (1+t)^{-\frac{p-1}{2 p}}.
\end{equation*}}
\noindent In particular, $ \|U^r(t,x)-u^r(\frac{x}{t})\|_{L^\infty{}} \leq  C_{p, \delta} (1+t)^{-\frac{1}{2}}$.
\end{lemma}

\noindent The proof of (1)-(4) in Lemma \ref{Properties of the smooth rarefaction wave 1} could be found in \cite{wang2002stability}, and please refer \cite{zhao2000nonlinear} for (5).  Meanwhile, we have the uniform energy estimate for the perturbation 
$$\phi(t, x):=u(t, x)-{U}^r(t, x), \psi(t, x):=v(t, x)-f({U^r}(t, x)),$$
\noindent with the initial data $(\phi_0, \psi_0)(x)$.

\begin{lemma}[Uniform energy estimate \textrm{I} \cite{zhao2000nonlinear}]
\label{Uniform energy estimate: Stri-Convex}
It holds that there exists a constant $C\left(\phi_0, \psi_0\right)>0$ depending on $\left\|\phi_0\right\|_{H^2},\left\|\psi_0\right\|_{H^2}$ and on the fixed flux and wave profiles, such that for each $0<t<+\infty$, the unique global smooth solution $\left(u(t, x), v(t, x)\right)$ in 
Lemma \ref{Global existence} satisfies
\begin{equation*} \|\phi(t)\|_{H^2}^2+\left\|\phi_t(t)\right\|_{H^1}^2+\int_0^t\left\|\phi_x(\tau)\right\|_{H^1}^2+\left\|\phi_t(\tau)\right\|_{H^1}^2\mathrm{d} \tau \leq C\left(\phi_0, \psi_0\right), 
\end{equation*}
\noindent and
\begin{equation*} \|\psi(t)\|_{H^2}^2+\left\|\psi_t(t)\right\|_{H^1}^2+\int_0^t\left\|\psi_x(\tau)\right\|_{H^1}^2+\left\|\psi_t(\tau)\right\|_{H^1}^2\mathrm{d} \tau \leq C\left(\phi_0, \psi_0\right) .
\end{equation*}
\end{lemma}

We now restate one of the main results.
\begin{theorem}[Strictly convex case]
\label{Main Th 1}
Suppose that the flux $f\in C^4(\mathbb{R})$ is strictly convex and satisfies the sub-characteristic condition (\ref{sub-characteristic}). If $u_0-U_0^r, v_0-f\left(U_0^r\right) \in L^1(\mathbb{R}) \cap H^2(\mathbb{R})$, then the Cauchy problem
(\ref{Original equation: Jin-Xin system})-(\ref{initial values}) admits a unique global smooth solution $(u(t,x), v(t,x))$,  which satisfies
\begin{equation*}
\sup_{t\geq0}\|u(t)\|_{L^\infty}\leq2N_0,\qquad \sup_{t\geq0}\|v(t)\|_{L^\infty}\leq B(N_0).
\end{equation*}
\noindent Here, the state space $M$ in \eqref{sub-characteristic} is $[-2N_0,2N_0]$, where
$$
N_0=\max \left\{\left\|u_0(x)\right\|_{L^{\infty}},\left\|v_0(x)\right\|_{L^{\infty}{}}\right\}, F\left(N_0\right)=\sup _{|u| \leq N_0}\left|f(u)\right|, B\left(N_0\right)=2 N_0+F\left(2 N_0\right) .
$$

Moreover, the solution satisfies the following properties:
\begin{enumerate}[label=(\roman*), leftmargin=*]
    \item ({\(L^1\) bound}). There exists a constant $C>0$, 
such that for $t>0$,
\begin{equation*}
    \left\|\left(u(t, x)-U^r(t,x), v(t, x)-f\left(U^r(t,x)\right)\right)\right\|_{L^1{}} \leq C\bigl[1+\ln(1+t)\bigr].
\end{equation*}

    \item ({Basic \(L^2\) decay}). For any fixed $\delta>0$, there exists a constant $C>0$, such that for $t>0$,
\begin{equation*}
 \left\|\left(u-U^r, v-f\left(U^r\right)\right)(t, x)\right\|_{H^1} \leq C(1+t)^{-\frac{1}{4}+\frac{\delta}{2}}.
\end{equation*}

    \item ({\(L^\infty\) decay}). For fixed $\delta>0$ as in (ii), there exists a constant $C>0$, such that for $t>0$,
\begin{equation*}
  \left\|\left(u(t, x)-u^r\left(\frac{x}{t}\right), v(t, x)-f\left(u^r\left(\frac{x}{t}\right)\right)\right)\right\|_{L^{\infty}{}} \leq C(1+t)^{-\frac{1}{2}+\delta},
\end{equation*}
\end{enumerate}
\end{theorem}
\begin{remark}
Liu \cite{liu2001lp} first established the basic \(L^2\) decay estimate (ii) via the time-weighted \(L^2\) method, assuming the \(L^1\) contraction estimate (i). Building on Liu's work, we employ a refined decomposition method to verify his conjecture on the \(L^1\) bound of the perturbation. Furthermore, after converting \eqref{Original equation: Jin-Xin system} into an equivalent damped wave equation, we derive the basic \(L^2\) decay estimate for \((u-U^r, v-f(U^r))\) and the derivative \(L^2\) decay estimate for \(((u-U^r)_x, (v-f(U^r))_x)\). Consequently, the \(L^\infty\) decay rate follows by applying the Sobolev inequality.
\end{remark}

\subsection{The $L^1$ estimate of $(\phi, \psi)$}
Next, we introduce the $L^1$ estimate, where we use a decomposition method and the standard $L^1$ estimate approach to obtain logarithmic control motivated by Liu\cite{liu2001lp}. The corresponding \(L^1\) estimate, which constitutes one of the main technical contributions of this work, is as follows:
\begin{lemma}
\label{L1 estimate: Stri-convex}
Under the same assumptions as in Theorem \ref{Main Th 1}, there exists a constant $C>0$, 
such that for $t\geq0$,
\begin{equation*}
    \left\|\left(\phi(t, x), \psi(t, x)\right)\right\|_{L^1{}} \leq C\bigl[1+\ln(1+t)\bigr]\leq C_\delta (1+t)^\delta,
\end{equation*}
\noindent where \(\delta>0\) is a fixed constant and \(C_\delta>0\) is a constant depending on \(\delta\).
\end{lemma}
\begin{proof}
At first, we convert the equation that the smooth rarefaction wave $U^r(t,x)$ satisfies to the approximate Jin-Xin relaxation system. Here, compared to the common estimate, we apply a decomposition method.  Denote $V^r := f(U^r) - \left(1 - (f'(U^r))^2\right) U^r_x$, we have 
\begin{equation}
    U_t^r+V_x^r=-\left(\left(1-\left(f^{\prime}\left(U^r\right)\right)^2\right) U_x^r\right)_x.
\end{equation}
\noindent By direct computation, we obtain that
\begin{equation}
\label{Jin-Xin L^1 eq1.2}
    \begin{aligned}
& V_t^r+U_x^r-\left(f\left(U^r\right)-V^r\right) \\
= & f^{\prime}\left(U^r\right) U_t^r-\left(\left(1-\left(f^{\prime}\left(U^r\right)\right)^2\right) U_x^r\right)_t+U_x^r-\left(1-\left(f^{\prime}\left(U^r\right)\right)^2\right) U_x^r \\
= & f^{\prime}\left(U^r\right)\left(-f^{\prime}\left(U^r\right) U_x^r\right)-\left(\left(1-\left(f^{\prime}\left(U^r\right)\right)^2\right) U_x^r\right)_t+U_x^r-\left(1-\left(f^{\prime}\left(U^r\right)\right)^2\right) U_x^r \\
= & -\left(\left(1-\left(f^{\prime}\left(U^r\right)\right)^2\right) U_x^r\right)_t.
\end{aligned}
\end{equation}
\noindent Let the residuals $e_1:=-\left(\left(1-\left(f^{\prime}\left(U^r\right)\right)^2\right) U_x^r\right)_x$ and $e_2:=-\left(\left(1-\left(f^{\prime}\left(U^r\right)\right)^2\right) U_x^r\right)_t$, hence, (\ref{Jin-Xin L^1 eq1.2}) give that 
\begin{equation}
\label{Jin-Xin L^1 eq1.3}
    \left\{\begin{array}{l}
U_t^r+V_x^r=e_1, \\
V_t^r+U_x^r-\left(f\left(U^r\right)-V^r\right)=e_2 .
\end{array}\right.
\end{equation}
\noindent It is obvious that (\ref{Jin-Xin L^1 eq1.3}) is the Jin-Xin model with residuals for $(U^r, V^r)$.  Next, by the standard Riemann invariants estimate procedure, let 
$$\varPsi:=v-V^r=\psi+\left(1 - (f'(U^r))^2\right) U^r_x.$$
\noindent Due to $U^r_x(t,x)\in L^\infty(L^1(\mathbb{R}))$, $\psi_{0}(x):=\psi|_{t=0}\in L^1(\mathbb{R})$ leads to $\varPsi_{0}(x):=\varPsi|_{t=0}\in L^1(\mathbb{R})$.
It holds that from (\ref{Original equation: Jin-Xin system}) and (\ref{Jin-Xin L^1 eq1.3})
\begin{equation}
\label{Jin-Xin L^1 perturbation system}
    \left\{\begin{array}{l}
\phi_t+\varPsi_{ x}=-e_1, \\
\varPsi_t+\phi_x=f\left(U^r+\phi\right)-f\left(U^r\right)-\varPsi-e_2 .
\end{array}\right.
\end{equation}
\noindent Let $R:=\varPsi-\phi$, $S:=-\varPsi-\phi$, so we can diagonalize (\ref{Jin-Xin L^1 perturbation system}) to 
\begin{equation}
    \label{Jin-Xin L^1 eq1.4}
    \left\{\begin{array}{l}
R_t-R_x=-\frac{R}{2}+\frac{S}{2}+\left(f\left(U^r+\phi\right)-f\left(U^r\right)\right)-e_2+e_1, \\
S_t+S_x=-\frac{S}{2}+\frac{R}{2}-\left(f\left(U^r+\phi\right)-f\left(U^r\right)\right)+e_2+e_1.
\end{array}\right.
\end{equation}
\noindent The mean-value formula tells us that 
\begin{equation}
    \begin{aligned}
f\left(U^r+\phi\right)-f\left(U^r\right) & =f^{\prime}(\eta) {} \phi
 =-\frac{f^{\prime}(\eta)}{2} R-\frac{f^{\prime}(\eta)}{2} S,
\end{aligned}
\end{equation}
\noindent where $\eta$ is between $U^r$ and $U^r+\phi$. Therefore,
\begin{equation}
\label{Jin-Xin L^1 eq1.5}
    \left\{\begin{array}{l}
R_t-R_x=-\frac{1+f^{\prime}(\eta)}{2} R+\frac{1-f^{\prime}(\eta)}{2} S-e_2+e_1, \\
S_t+S_x=-\frac{1-f^{\prime}(\eta)}{2} S+\frac{1+f^{\prime}(\eta)}{2} R+e_2+e_1.
\end{array}\right.
\end{equation}
\noindent Now multiply $(\ref{Jin-Xin L^1 eq1.5})_{1}$ by $\operatorname{sgn} R$ and $(\ref{Jin-Xin L^1 eq1.5})_{2}$ by $\operatorname{sgn} S$, then sum them up, one could know
\begin{equation}
    \begin{aligned}
& (|R|)_t+(|S|)_t+\{-|R|+|S|\}_x \\
= & -(\operatorname{sgn} R-\operatorname{sgn} S)\left(\frac{1+f^{\prime}(\eta)}{2} R-\frac{1-f^{\prime}(\eta)}{2} S\right) \\
& +\left(e_1-e_2\right) \operatorname{sgn} R+\left(e_1+e_2\right) \operatorname{sgn} S.
\end{aligned}
\end{equation}
\noindent The sub-characteristic condition (\ref{sub-characteristic}) and the definition of the sgn function together imply that 
\begin{equation}
    (\operatorname{sgn} R-\operatorname{sgn} S)\left(\frac{1+f^{\prime}(\eta)}{2} R-\frac{1-f^{\prime}(\eta)}{2} S\right)\geq0.
\end{equation}
\noindent Therefore, one could have
\begin{equation}
\label{Jin-Xin L^1 eq1.6}
    (|R|)_t+(|S|)_t+\{-|R|+|S|\}_x\leq\left(e_1-e_2\right) \operatorname{sgn} R+\left(e_1+e_2\right) \operatorname{sgn} S.
\end{equation}
\noindent Integrating (\ref{Jin-Xin L^1 eq1.6}) over $\mathbb{R}$, we arrive at
\begin{equation}
\|R(t)\|_{L^1}+\|S(t)\|_{L^1} \leq\|R(0)\|_{L^1}+\|S(0)\|_{L^1}+2 \int_0^t\left(\left\|e_1(\tau)\right\|_{L^1}+\left\|e_2(\tau)\right\|_{L^1}\right) d \tau.
\end{equation}
\noindent From Lemma \ref{Properties of the smooth rarefaction wave 1}, we have $\left\|e_1( t)\right\|_{L^1}, \left\|e_2( t)\right\|_{L^1}\leq C(1+t)^{-1}$. Thus, we can deduce that 
\begin{equation}
    \begin{aligned}
\|R(t)\|_{L^1}+\|S(t)\|_{L^1} \leq & \|R(0)\|_{L^1}+\|S(0)\|_{L^1}+C \int_0^t(1+\tau)^{-1} d \tau \\
& \leq C\bigl[1+\ln(1+t)\bigr],
\end{aligned}
\end{equation}
\noindent that is 
\begin{equation}
    \|\phi(t)\|_{L^1}+\|\varPsi(t)\|_{L^1}\leq C\bigl[1+\ln(1+t)\bigr].
\end{equation}
\noindent By $\|U^r_x(t)\|_{L^1}\leq C<+\infty$, 
\begin{equation}
    \|\phi(t)\|_{L^1}+\|\psi(t)\|_{L^1}\leq C\bigl[1+\ln(1+t)\bigr].
\end{equation}
\noindent  Consequently, we complete the proof of the $L^1$ estimate.
The estimate here is relative to the smooth rarefaction wave of Type 1. 
\end{proof}

\subsection{$L^2$ decay for $(\phi, \psi)$}
In the following parts, we only need to get the decay rate of $\phi$. Indeed, according to (\ref{Jin-Xin L^1 perturbation system}), we have
\begin{equation*}
    \begin{aligned}
\varPsi_t & =-\varPsi+f\left(U^r+\phi\right)-f\left(U^r\right)-\phi_x-e_2 \\
& =-\varPsi+\left(O(1) \phi-\phi_x-e_2\right),
\end{aligned}
\end{equation*}
\noindent which leads to
\begin{equation}
\label{expression of psi_1}
    \varPsi(t, x)=e^{-t} \varPsi_0(x)+\int_0^t\left(O(1) \phi-\phi_x-e_2\right)(\tau, x) e^{-(t-\tau)} d \tau.
\end{equation}
\noindent Thus one could obtain the decay of $\psi$ by the fact that $\psi=\varPsi-O(1) U_x^r$. Now we derive the decay rate of $\|\phi(t)\|_{L^2}$ by applying the time-weighted $L^2$ method. 
\begin{lemma}
    \label{L2 decay : stri-convex}
Under the same assumptions as in Theorem \ref{Main Th 1}, for $\delta>0$ as in Lemma \ref{L1 estimate: Stri-convex} and $\gamma>2$,
there exists a constant $C>0$, 
such that for $t\geq0$,
\begin{equation*}
 \begin{aligned}
& (1+t)^\gamma\left(\|\phi(t)\|_{L^2}^2+\left\|\phi_t(t)\right\|_{L^2}^2+\left\|\phi_x(t)\right\|_{L^2}^2\right)+\int_0^t(1+\tau)^\gamma\left(\left\|\phi_t(\tau)\right\|_{L^2}^2+\left\|\phi_x(\tau)\right\|_{L^2}^2+\left\|\phi \sqrt{U_x^r}(\tau)\right\|_{L^2}^2\right) d \tau \\
& \leq C\left(\left\|\phi_0\right\|_{L^2}^2+\left\|\phi_{0 x}\right\|_{L^2}^2+\left\|\phi_1\right\|_{L^2}^2\right)+C(1+t)^{\gamma+\delta-\frac{1}{2}}.
\end{aligned}
\end{equation*}
\end{lemma}

\begin{proof}
First, eliminating $v$ from (\ref{Original equation: Jin-Xin system}) gives the scalar damped wave equation with nonlinear convection term
 \begin{equation}
 \label{damped wave equ}
     u_{t t}+u_t+f(u)_x=u_{x x}.
 \end{equation}

\noindent Next, reformulate the problem (\ref{damped wave equ}) to the perturbation equation according to the difference between the solution \(u(t,x)\) to \eqref{Original equation: Jin-Xin system} and the smooth 
rarefaction wave  $U^r(t, x)=\left(f^{\prime}\right)^{-1}(W(t, x))$, here
the smooth approximation $W(t, x)$ satisfies (\ref{2}). According to Lemma \ref{Properties of the smooth rarefaction wave 1} 
we can rewrite the problem \eqref{damped wave equ} in terms of $\phi$ as
\begin{equation}
\label{pertur-equ: stri-convex}
\begin{cases}
\phi_{tt} + \phi_t + \left( f(U^r+\phi) - f(U^r) \right)_x - \phi_{xx} = F(U^r), \\
\phi|_{t=0} = u_0(x) - U^r_0(x) =: \phi_0(x), \\
\phi_t|_{t=0} = u_1(x) - U^r_{0t}(x) =: \phi_1(x).
\end{cases}
\end{equation}
\noindent Where the residual 
$F(U^r):=U^r_{xx}-U^r_{tt}$. Here $u_1=-v_{0x}$ follows from the first equation of system (\ref{Original equation: Jin-Xin system}). Now multiplying (\ref{pertur-equ: stri-convex}) by $\phi+2\phi_t$,
\begin{equation}
\label{div-form1:stri-convex}
    \begin{aligned}
& \left\{\frac{1}{2} \phi^2+\phi \phi_t+\phi_t^2+\phi_x^2\right\}_t+\phi_t^2+\phi_x^2+\left\{-\phi \phi_x-2 \phi_t \phi_x\right\}_x \\
& +\phi {}\left(f\left(U^r+\phi\right)-f\left(U^r\right)\right)_x+2 \phi_t {}\left(f\left(U^r+\phi\right)-f\left(U^r\right)\right)_x=\phi {} F+2 \phi_t {} F .
\end{aligned}
\end{equation}
\noindent According to $2 \phi_t {}\left(f\left(U^r+\phi\right)-f\left(U^r\right)\right)_x=2 \phi_t {}\left(f^{\prime}\left(U^r+\phi\right)-f^{\prime}\left(U^r\right)\right) {U}^r_x+2 \phi_t {} f^{\prime}\left(U^r+\phi\right) {} \phi_x$, thus we can get
\begin{equation}
    \begin{aligned}
& \left\{\frac{1}{2} \phi^2+\phi \phi_t+\phi_t^2+\phi_x^2\right\}_t+\left(\phi_t^2+\phi_x^2+2 f^{\prime}(u) \phi_t \phi_x\right)+\left\{-\phi \phi_x-2 \phi_t \phi_x\right\}_x \\
& +\phi {}\left(f\left(U^r+\phi\right)-f\left(U^r\right)\right)_x+2 \phi_t\left(f^{\prime}\left(U^r+\phi\right)-f^{\prime}\left(U^r\right)\right) {U}^r_x=\phi {} F+2 \phi_t {} F .
\end{aligned}
\end{equation}
\noindent Notice that  $\phi_t^2+\phi_x^2+2 f^{\prime}(u) \phi_t \phi_x$ is a quadratic form with respect to $\phi_x$ and $\phi_t$, so we just need to prove its positivity. Since
$\phi_t^2+\phi_x^2+2 f^{\prime}(u) \phi_t \phi_x \geq \phi_t^2+\phi_x^2-2 k\left|\phi_t\right| {}\left|\phi_x\right|$, 
denote the matrix corresponding to the quadratic form
as 
$A=\left(\begin{array}{cc}
1 & -k \\
-k & 1
\end{array}\right)$. Thus $\mathrm{det} A=1-k^2>0$, where we have used the sub-characteristic condition (\ref{sub-characteristic}). Then there exists $\alpha>0$, such that 
\begin{equation}
    \label{positive-quard: stri-convex}
     \phi_t^2+\phi_x^2+2 f^{\prime}(u) \phi_t \phi_x\geq\alpha(\phi_t^2+\phi_x^2).
\end{equation}
\noindent Moreover, by using the strict convexity of \(f\),
\begin{equation}
    \label{dissipation: stri-convex}
   \begin{aligned}
& \phi {}\left(f\left(U^r+\phi\right)-f\left(U^r\right)\right)_x \\
= & \left\{\phi {}\left(f\left(U^r+\phi\right)-f\left(U^r\right)\right)-\int_0^\phi\left(f\left(U^r+\eta\right)-f\left(U^r\right)\right) d \eta\right\}_x+\left(f\left(U^r+\phi\right)-f\left(U^r\right)-f^{\prime}\left(U^r\right) {} \phi\right) U_x^r \\
\geq & \left\{\phi {}\left(f\left(U^r+\phi\right)-f\left(U^r\right)\right)-\int_0^\phi\left(f\left(U^r+\eta\right)-f\left(U^r\right)\right) d \eta\right\}_x+\frac{\kappa_0}{2} \phi^2 U_x^r.
\end{aligned}
\end{equation}

\noindent Substituting (\ref{positive-quard: stri-convex})  and (\ref{dissipation: stri-convex}) into (\ref{div-form1:stri-convex}), we deduce that 
\begin{equation}
\label{div-form2: stri convex}
    \begin{gathered}
\left\{\frac{1}{2} \phi^2+\phi \phi_t+\phi_t^2+\phi_x^2\right\}_t+\left\{\phi {}\left(f\left(U^r+\phi\right)-f\left(U^r\right)\right)-\int_0^\phi\left(f\left(U^r+\eta\right)-f\left(U^r\right)\right) d \eta-\phi \phi_x-2 \phi_t \phi_x\right\}_x \\
+\frac{\kappa_0}{2} \phi^2 U_x^r+\alpha\left(\phi_t^2+\phi_x^2\right)+2 \phi_t\left(f^{\prime}\left(U^r+\phi\right)-f^{\prime}\left(U^r\right)\right) U_x^r \leq \phi F+2 \phi_t F.
\end{gathered}
\end{equation}
\noindent Multiplying (\ref{div-form2: stri convex})   by $(1+t)^\gamma$, where $\gamma>0$ large enough is to be chosen, and then integrating it over $[0, t]\times \mathbb{R}$, one could get
\begin{align}
\label{time-weighted initial: stri convex}
& (1+t)^\gamma\left\|\left(\phi, \phi_x, \phi_t\right)(t)\right\|_{L^2}^2+\int_0^t(1+\tau)^\gamma\left(\left\|\phi_x(\tau)\right\|_{L^2}^2+\left\|\phi_t(\tau)\right\|_{L^2}^2+\left\|\phi \sqrt{U_x^r}(\tau)\right\|_{L^2}^2\right) d \tau \notag\\
\leq & C\left(\left\|\phi_0\right\|_{H^1}^2+\left\|\phi_1\right\|_{L^2}^2\right)+C \int_0^t(1+\tau)^{\gamma-1}\|\phi(\tau)\|_{L^2}^2 d \tau+C \int_0^t(1+\tau)^{\gamma-1}\left(\left\|\phi_t(\tau)\right\|_{L^2}^2+\left\|\phi_x(\tau)\right\|_{L^2}^2\right) d \tau \notag\\
& +C \int_0^t(1+\tau)^\gamma \int_{\mathbb{R}}\left|\phi_t\left(f^{\prime}\left(U^r+\phi\right)-f^{\prime}\left(U^r\right)\right) U_x^r\right|(\tau, x) d x d \tau+C \int_0^t(1+\tau)^\gamma \int_{\mathbb{R}}|\phi||F| d x d \tau \notag\\
& +C \int_0^t(1+\tau)^\gamma \int_{\mathbb{R}}\left|\phi_t\right||F| d x d \tau=: C\left(\left\|\phi_0\right\|_{H^1}^2+\left\|\phi_1\right\|_{L^2}^2\right)+\sum_{i=1}^5 I_i.
\end{align}
Now, we estimate $I_1-I_5$ term by term. First, 
 applying Lemma \ref{L1 estimate: Stri-convex} with exponent $\delta/2$, we have $\|\phi(t)\|_{L^1}\leq C_\delta(1+t)^{\delta/2}$. Thus the Nash inequality  gives
\begin{equation}
    \|\phi(t)\|_{L^2} \leq C\left\|\phi_x(t)\right\|{ }_{L^2}^{1 / 3}\|\phi(t)\|_{L^1}^{2 / 3} \leq C\left\|\phi_x(t)\right\|{ }_{L^2}^{1 / 3} {}(1+t)^{\delta/3},
\end{equation}
\noindent which leads to
    \begin{align}
    \label{I1 esti}
I_1 & \leq C \int_0^t(1+\tau)^{\gamma-1+\frac{2}{3} \delta} {}\left\|\phi_x(\tau)\right\|_{L^2}^{2 / 3} d \tau \notag\\
& =C \int_0^t(1+\tau)^{\frac{2 \gamma}{3}-1+\frac{2 \delta}{3}} {}\left((1+\tau)^\gamma\left\|\phi_x(\tau)\right\|_{L^2}^2\right)^{1 / 3} d \tau \notag\\
& \leq \varepsilon \int_0^t(1+\tau)^\gamma\left\|\phi_x(\tau)\right\|_{L^2}^2 d \tau+C \int_0^t(1+\tau)^{\gamma+\delta-\frac{3}{2}} d \tau \notag\\
& \leq \varepsilon \int_0^t(1+\tau)^\gamma\left\|\phi_x(\tau)\right\|_{L^2}^2 d \tau+C(1+t)^{\gamma+\delta-\frac{1}{2}},
\end{align}
\noindent where we have used the Young inequality with $\frac{1}{3}+\frac{2}{3}=1$ and taken $\gamma>3/2$. Next, owing to Lemma \ref{Uniform energy estimate: Stri-Convex}, take $\gamma>1$, we deduce that 
\begin{equation}
\label{I2 esti}
    I_2 \leq C(1+t)^{\gamma-1} \int_0^t\left(\left\|\phi_t\right\|_{L^2}^2+\left\|\phi_x\right\|_{L^2}^2\right) d \tau \leq C(1+t)^{\gamma-1}.
\end{equation}
\noindent Subsequently, take $\gamma>2$, the Young inequality, Lemma \ref{Properties of the smooth rarefaction wave 1} and Lemma \ref{Uniform energy estimate: Stri-Convex} together imply that 
\begin{align}
\label{I3 esti}
I_3 & \leq C \int_0^t(1+\tau)^\gamma \int_{\mathbb{R}}\left|\phi_t\right| {}\left(|\phi| {} U_x^r\right) d x d \tau \notag\\
& \leq \varepsilon \int_0^t(1+\tau)^\gamma\left\|\phi_t\right\|_{L^2}^2 d \tau+C \int_0^t(1+\tau)^\gamma \int_{\mathbb{R}} \phi^2 {}\left(U_x^r\right)^2 d x d \tau \notag\\
& \leq \varepsilon \int_0^t(1+\tau)^\gamma\left\|\phi_t\right\|_{L^2}^2 d \tau+C \int_0^t(1+\tau)^\gamma\left\|U_x^r(\tau)\right\|_{L^{\infty}}^2\|\phi\|_{L^2}^2 d \tau \notag\\
& \leq \varepsilon \int_0^t(1+\tau)^\gamma\left\|\phi_t\right\|_{L^2}^2 d \tau+C \int_0^t(1+\tau)^{\gamma-2}\|\phi\|_{L^2}^2 d \tau \notag\\
& \leq \varepsilon \int_0^t(1+\tau)^\gamma\left\|\phi_t\right\|_{L^2}^2 d \tau+C \int_0^t(1+\tau)^{\gamma-2} d \tau \notag\\
& \leq \varepsilon \int_0^t(1+\tau)^\gamma\left\|\phi_t\right\|_{L^2}^2 d \tau+C(1+t)^{\gamma-1}.  
\end{align}
\noindent As to $I_4$, Lemma \ref{Properties of the smooth rarefaction wave 1} tells us that 
\begin{equation}
    \|F(U^r)(t)\|_{L^1} \leq C(\left\|U^r_{tt}\right\|_{L^1}+\left\|U^r_{x x}\right\|_{L^1}) \leq C(1+t)^{-1},
\end{equation}
\noindent thus,  by using the interpolation inequality,
    \begin{align}
        \label{I4 esti}
I_4 & \leq C \int_0^t(1+\tau)^\gamma \int_{\mathbb{R}}|\phi| {}|F| d x d \tau \leq C \int_0^t(1+\tau)^\gamma\|\phi(\tau)\|_{L^{\infty}} {}\|F(\tau)\|_{L^1} d \tau \notag\\
& \leq C \int_0^t(1+\tau)^\gamma\|\phi(\tau)\|_{L^1}^{1 / 3} {}\left\|\phi_x(\tau)\right\|_{L^2}^{2 / 3} {}\|F(\tau)\|_{L^1} d \tau \notag\\
& \leq C \int_0^t(1+\tau)^{\gamma-1+\frac{1}{3} \delta}\left\|\phi_x(\tau)\right\|_{L^2}^{2 / 3} d \tau \notag\\
& \leq \varepsilon \int_0^t(1+\tau)^\gamma\left\|\phi_x(\tau)\right\|_{L^2}^2 d \tau+C(1+t)^{\gamma+\delta-\frac{1}{2}},
\end{align}
\noindent here we have used the same argument as in the estimate of (\ref{I1 esti}). Finally, from the Young inequality and Lemma \ref{Properties of the smooth rarefaction wave 1},
\begin{equation}
\label{I5 esti}
    \begin{aligned}
I_5 & \leq \varepsilon \int_0^t(1+\tau)^\gamma\left\|\phi_t(\tau)\right\|_{L^2}^2 d \tau+C \int_0^t(1+\tau)^\gamma\|F\|{ }_{L^2}^2 d \tau \\
& \leq \varepsilon \int_0^t(1+\tau)^\gamma\left\|\phi_t(\tau)\right\|_{L^2}^2 d \tau+C \int_0^t(1+\tau)^{\gamma-\frac{5}{2}} d \tau \\
& \leq \varepsilon \int_0^t(1+\tau)^\gamma\left\|\phi_t(\tau)\right\|_{L^2}^2 d \tau+C \int_0^t(1+\tau)^{\gamma-\frac{3}{2}} d \tau \\
& \leq \varepsilon \int_0^t(1+\tau)^\gamma\left\|\phi_t(\tau)\right\|_{L^2}^2 d \tau+C(1+t)^{\gamma-\frac{1}{2}}.
\end{aligned}
\end{equation}
Consequently, inserting (\ref{I1 esti}), (\ref{I2 esti}), (\ref{I3 esti}), (\ref{I4 esti}) and (\ref{I5 esti}) into (\ref{time-weighted initial: stri convex}) and then taking $\varepsilon>0$ small enough, it can be obtained that for any fixed $\delta>0$, $\gamma>2$, 
\begin{equation}
 \begin{aligned}
& (1+t)^\gamma\left(\|\phi(t)\|_{L^2}^2+\left\|\phi_t(t)\right\|_{L^2}^2+\left\|\phi_x(t)\right\|_{L^2}^2\right)+\int_0^t(1+\tau)^\gamma\left(\left\|\phi_t(\tau)\right\|_{L^2}^2+\left\|\phi_x(\tau)\right\|_{L^2}^2+\left\|\phi \sqrt{U_x^r}(\tau)\right\|_{L^2}^2\right) d \tau \\
& \leq C\left(\left\|\phi_0\right\|_{L^2}^2+\left\|\phi_{0 x}\right\|_{L^2}^2+\left\|\phi_1\right\|_{L^2}^2\right)+C(1+t)^{\gamma+\delta-\frac{1}{2}}.
\end{aligned}
\end{equation}
\end{proof}

\subsection{$L^2$ decay for $(\phi_x, \psi_x)$}
In Lemma \ref{L2 decay : stri-convex}, we have obtained the decay of $\left\|\phi(t)\right\|_{L^2}$ and $\left\|\phi_x(t)\right\|_{L^2}$. However, we can give a further estimate on $\left\|\phi_x(t)\right\|_{L^2}$. Next, we start to give the precise $L^2$ decay for $\|\phi_x(t)\|_{L^2}$. In other words, we improve the decay results of Liu in \cite{liu2001lp}.
\begin{lemma}
    \label{derivative L2 decay: stri-convex}
Under the same assumptions as in Theorem \ref{Main Th 1}, for $\delta>0$  as in Lemma \ref{L1 estimate: Stri-convex} and $\gamma>5/2$,
there exists a constant $C>0$, 
such that for $t\geq0$,
\begin{equation*}
\begin{aligned}
& (1+t)^\gamma\left(\left\|\phi_x(t)\right\|_{L^2}^2+\left\|\phi_{t x}(t)\right\|_{L^2}^2+\left\|\phi_{x x}(t)\right\|_{L^2}^2\right)+\int_0^t(1+\tau)^\gamma\left(\left\|\phi_{t x}(\tau)\right\|_{L^2}^2+\left\|\phi_{x x}(\tau)\right\|_{L^2}^2\right) d \tau \\
& \leq C\left(\left\|\phi_{0 x}\right\|_{L^2}^2+\left\|\phi_{0 x x}\right\|_{L^2}^2+\left\|\phi_{1 x}\right\|_{L^2}^2\right)+C(1+t)^{\gamma-\frac{3}{2}+5\delta}.
\end{aligned}.
\end{equation*}
\end{lemma}
\begin{proof}
Taking weak derivative on  (\ref{pertur-equ: stri-convex}) with respect to $x$, we have
\begin{equation}
\label{pertu equ derivative: stri convex}
\phi_{t t x}+\phi_{t x}+(f(U^r+\phi)-f(U^r))_{x x}-\phi_{x x x}=F(U^r)_x.
\end{equation}
\noindent Multiplying (\ref{pertu equ derivative: stri convex}) by $\phi_x+2\phi_{tx}$, it follows that 
\begin{equation}
    \begin{aligned}
\left\{\phi_x \phi_{t x}+\frac{1}{2} \phi_x^2+\phi_{t x}^2+\phi_{x x}^2\right\}_t & +\phi_x\left(f\left(U^r+\phi\right)-f\left(U^r\right)\right)_{x x}+2 \phi_{t x} {}\left(f\left(U^r+\phi\right)-f\left(U^r\right)\right)_{x x} \\
& +\phi_{t x}^2+\phi_{x x}^2+\left\{-\phi_x \phi_{x x}-2 \phi_{t x} \phi_{x x}\right\}_x=\phi_x F_x+2 \phi_{t x} F_x.
\end{aligned}
\end{equation}
\noindent A straightforward calculation implies that 
\begin{equation}
    \begin{aligned}
& 2 \phi_{t x} {}\left(f\left(U^r+\phi\right)-f\left(U^r\right)\right)_{x x} \\
= & 2 f^{\prime}(u) \phi_{t x} {} \phi_{x x}+2 \phi_{t x} {}\left(f^{\prime}\left(U^r+\phi\right)-f^{\prime}\left(U^r\right)\right) U_{x x}^r+2 \phi_{t x} {}\left(f^{\prime \prime}\left(U^r+\phi\right)-f^{\prime \prime}\left(U^r\right)\right)\left(U_x^r\right)^2 \\
& +2 \phi_{t x} f^{\prime \prime}(u) {} \phi_x^2+4 \phi_{t x} {} f^{\prime \prime}(u){} \phi_x {} U_x^r,
\end{aligned}
\end{equation}
\noindent and
\begin{equation}
    \begin{aligned}
& \quad\phi_x {}\left(f\left(U^r+\phi\right)-f\left(U^r\right)\right)_{x x} \\
&=\left\{\phi_x {}\left(f\left(U^r+\phi\right)-f\left(U^r\right)\right)_x\right\}_x-\phi_{x x} {}\left(f^{\prime}\left(U^r+\phi\right)-f^{\prime}\left(U^r\right)\right) U_x^r-\phi_{x x} {}\left(f^{\prime}\left(U^r+\phi\right)-f^{\prime}\left(U^r\right)\right) \phi_x \\
&\quad-\frac{1}{2}\left\{\phi_x^2 {} f^{\prime}\left(U^r\right)\right\}_x+\frac{1}{2} \phi_x^2 {} f^{\prime \prime}\left(U^r\right) U_x^r.
\end{aligned}
\end{equation}
\noindent Thus using the strict convexity of \(f\), one could have
\begin{align}
\label{derivative div form: stri-convex}
& \left\{\phi_x \phi_{t x}+\frac{1}{2} \phi_x^2+\phi_{t x}^2+\phi_{x x}^2\right\}_t+\left\{-\phi_x \phi_{x x}-2 \phi_{t x} \phi_{x x}-\phi_x F+\phi_x \left(f\left(U^r+\phi\right)-f\left(U^r\right)\right)_x-\frac{1}{2} \phi_x^2  f^{\prime}\left(U^r\right)\right\}_x \notag\\
& +\alpha \left(\phi_{t x}^2+\phi_{x x}^2\right)+2 \phi_{t x} \left(f^{\prime}\left(U^r+\phi\right)-f^{\prime}\left(U^r\right)\right) U_{x x}^r+2 \phi_{t x}\left(f^{\prime \prime}\left(U^r+\phi\right)-f^{\prime \prime}\left(U^r\right)\right)\left(U_x^r\right)^2 \notag\\
& +2 \phi_{t x} f^{\prime \prime}(u)\phi_x^2+4 \phi_{t x} f^{\prime \prime}(u) \phi_x U_x^r-\phi_{x x} \left(f^{\prime}\left(U^r+\phi\right)-f^{\prime}\left(U^r\right)\right) U_x^r-\phi_{x x} \left(f^{\prime}\left(U^r+\phi\right)-f^{\prime}\left(U^r\right)\right) \phi_x \notag\\
& \leq-\phi_{x x} F+2 \phi_{t x} F_x.
\end{align}
\noindent where we applied 
\begin{equation}
     \phi_{tx}^2+\phi_{xx}^2+2 f^{\prime}(u) \phi_{tx} \phi_{xx}\geq\alpha(\phi_{tx}^2+\phi_{xx}^2).
\end{equation}
\noindent Now multiplying (\ref{derivative div form: stri-convex})   by $(1+t)^\gamma$, where $\gamma>0$ large enough is to be determined, and then integrating it over $[0, t]\times \mathbb{R}$, we can obtain that
\begin{align}
\label{derivative esti-time weighted initial-stri convex}
& (1+t)^\gamma\left(\left\|\phi_x(t)\right\|_{L^2}^2+\left\|\phi_{t x}(t)\right\|_{L^2}^2+\left\|\phi_{x x}(t)\right\|_{L^2}^2\right)+\int_0^t(1+\tau)^\gamma\left(\left\|\phi_{x x}(\tau)\right\|_{L^2}^2+\left\|\phi_{t x}(\tau)\right\|_{L^2}^2\right) d \tau \notag\\
\leq & C\left(\left\|\phi_{0 x}\right\|_{L^2}^2+\left\|\phi_{0 x x}\right\|_{L^2}^2+\left\|\phi_{1 x}\right\|_{L^2}^2\right)+C \int_0^t(1+\tau)^{\gamma-1}\left(\left\|\phi_x(\tau)\right\|_{L^2}^2+\left\|\phi_{t x}(\tau)\right\|_{L^2}^2+\left\|\phi_{x x}(\tau)\right\|_{L^2}^2\right) d \tau \notag\\
& +C \int_0^t(1+\tau)^\gamma\left|\int_{\mathbb{R}} \phi_{t x}\left(f^{\prime}\left(U^r+\phi\right)-f^{\prime}\left(U^r\right)\right) U_{x x}^r d x\right| d \tau \notag\\
& +C \int_0^t(1+\tau)^\gamma\left|\int_{\mathbb{R}} \phi_{t x}\left(f^{\prime \prime}\left(U^r+\phi\right)-f^{\prime \prime}\left(U^r\right)\right)\left(U_x^r\right)^2 d x\right| d \tau \notag\\
& +C \int_0^t(1+\tau)^\gamma\left|\int_{\mathbb{R}} \phi_{t x} f^{\prime \prime}\left(U^r+\phi\right) {} \phi_x^2 d x\right| d \tau+C \int_0^t(1+\tau)^\gamma\left|\int_{\mathbb{R}} \phi_{t x} f^{\prime \prime}\left(U^r+\phi\right) \phi_x U_x^r d x\right| d \tau \notag\\
& +C \int_0^t(1+\tau)^\gamma\left|\int_{\mathbb{R}} \phi_{x x}\left(f^{\prime}\left(U^r+\phi\right)-f^{\prime}\left(U^r\right)\right) U_x^r d x\right| d \tau \notag\\
& +C \int_0^t(1+\tau)^\gamma\left|\int_{\mathbb{R}} \phi_{x x}\left(f^{\prime}\left(U^r+\phi\right)-f^{\prime}\left(U^r\right)\right) \phi_x d x\right| d \tau \notag\\
& +C \int_0^t(1+\tau)^\gamma \int_{\mathbb{R}}\left|\phi_{x x}\right| {}|F| d x d \tau+C \int_0^t(1+\tau)^\gamma \int_{\mathbb{R}}\left|\phi_{t x}\right| {}\left|F_x\right| d x d \tau \notag\\
& =: C\left(\left\|\phi_{0 x}\right\|_{L^2}^2+\left\|\phi_{0 x x}\right\|_{L^2}^2+\left\|\phi_{1 x}\right\|_{L^2}^2\right)+\sum_{i=6}^{14} I_i.
\end{align}
\noindent Next we start to estimate $I_6\sim I_{14}$ term by term. First, using Lemma \ref{L2 decay : stri-convex}, Lemma \ref{Interpolation inequality} and the Young inequality, one could have 
    \begin{align}
    \label{phi_x L^2 weighted estimate}
\int_0^t(1+\tau)^{\gamma-1}\left\|\phi_x(\tau)\right\|_{L^2}^2 d \tau & \leq \int_0^t(1+\tau)^{\gamma-1}\left\|\phi_{x x}(\tau)\right\|_{L^2}\|\phi(\tau)\|_{L^2} d \tau \notag\\
& \leq C \int_0^t(1+\tau)^{\gamma-\frac{5}{4}+\frac{\delta}{2}}\left\|\phi_{x x}(\tau)\right\|_{L^2} d \tau \notag\\
& =C \int_0^t\left\{(1+\tau)^{\frac{\gamma}{2}}\left\|\phi_{x x}(\tau)\right\|_{L^2}\right\} {}(1+\tau)^{\frac{\gamma}{2}-\frac{5}{4}+\frac{\delta}{2}} d \tau \notag\\
& \leq \varepsilon \int_0^t(1+\tau)^\gamma\left\|\phi_{x x}(\tau)\right\|_{L^2}^2 d \tau+C \int_0^t(1+\tau)^{\gamma-\frac{5}{2}+\delta} d \tau \notag\\
& \leq \varepsilon \int_0^t(1+\tau)^\gamma\left\|\phi_{x x}(\tau)\right\|_{L^2}^2 d \tau+C(1+t)^{\gamma-\frac{3}{2}+\delta},
\end{align}
\noindent where we take $\gamma>5/2$. And thanks to Lemma \ref{Uniform energy estimate: Stri-Convex}, using the Young inequality, 
it follows that 
    \begin{align}
& \int_0^t(1+\tau)^{\gamma-1}\left(\left\|\phi_{t x}(\tau)\right\|_{L^2}^2+\left\|\phi_{x x}(\tau)\right\|_{L^2}^2\right) d \tau \notag\\
\leq & C \int_0^t\left\{(1+\tau)^{\frac{\gamma}{2}}\left(\left\|\phi_{t x}\right\|_{L^2}^2+\left\|\phi_{x x}\right\|_{L^2}^2\right)^{\frac{1}{2}}\right\} {}\left\{(1+\tau)^{\frac{\gamma}{2}-1}\left(\left\|\phi_{t x}\right\|_{L^2}^2+\left\|\phi_{x x}\right\|_{L^2}^2\right)^{\frac{1}{2}}\right\} d \tau \notag\\
\leq & \varepsilon \int_0^t(1+\tau)^\gamma\left(\left\|\phi_{t x}\right\|_{L^2}^2+\left\|\phi_{x x}\right\|_{L^2}^2\right) d \tau+C \int_0^t(1+\tau)^{\gamma-2}\left(\left\|\phi_{t x}\right\|_{L^2}^2+\left\|\phi_{x x}\right\|_{L^2}^2\right) d \tau \notag\\
\leq & \varepsilon \int_0^t(1+\tau)^\gamma\left(\left\|\phi_{t x}\right\|_{L^2}^2+\left\|\phi_{x x}\right\|_{L^2}^2\right) d \tau+C(1+t)^{\gamma-2} \int_0^t\left\|\phi_{t x}\right\|_{L^2}^2+\left\|\phi_{x x}\right\|_{L^2}^2 d \tau \notag\\
\leq & \varepsilon \int_0^t(1+\tau)^\gamma\left(\left\|\phi_{t x}\right\|_{L^2}^2+\left\|\phi_{x x}\right\|_{L^2}^2\right) d \tau+C(1+t)^{\gamma-2},
\end{align}
\noindent hence
\begin{equation}
    \label{I6 esti}
    I_6 \leq \varepsilon \int_0^t(1+\tau)^\gamma\left(\left\|\phi_{t x}\right\|_{L^2}^2+\left\|\phi_{x x}\right\|_{L^2}^2\right) d \tau+C(1+t)^{\gamma-\frac{3}{2}+\delta}.
\end{equation}
\noindent Then, Lemma \ref{Properties of the smooth rarefaction wave 1}, the mean-value formula and the Young inequality together lead to 
\begin{align}
  \label{I7, I8 esti}
I_7+I_8 \leq & \varepsilon \int_0^t(1+\tau)^\gamma\left(\left\|\phi_{t x}\right\|_{L^2}^2+\left\|\phi_{x x}\right\|_{L^2}^2\right) d \tau+C \int_0^t(1+\tau)^\gamma \int_{\mathbb{R}} \phi^2 {}\left(U_{x x}^r\right)^2 d x d \tau \notag\\
& +C \int_0^t(1+\tau)^\gamma \int_{\mathbb{R}} \phi^2 {}\left(U_x^r\right)^4 d x d \tau \notag\\
\leq & \varepsilon \int_0^t(1+\tau)^\gamma\left(\left\|\phi_{t x}\right\|_{L^2}^2+\left\|\phi_{x x}\right\|_{L^2}^2\right) d \tau+C \int_0^t(1+\tau)^\gamma\left(\left\|U_{x x}^r\right\|_{L^{\infty}}^2+\left\|U_x^r\right\|_{L^{\infty}}^4\right)\|\phi(\tau)\|_{L^2}^2 d \tau \notag\\
\leq & \varepsilon \int_0^t(1+\tau)^\gamma\left(\left\|\phi_{t x}\right\|_{L^2}^2+\left\|\phi_{x x}\right\|_{L^2}^2\right) d \tau+C \int_0^t(1+\tau)^{\gamma-3} d \tau \notag\\
\leq & \varepsilon \int_0^t(1+\tau)^\gamma\left(\left\|\phi_{t x}\right\|_{L^2}^2+\left\|\phi_{x x}\right\|_{L^2}^2\right) d \tau+C(1+t)^{\gamma-2}.
\end{align}
\noindent As to $I_9$, by applying Lemma \ref{Interpolation inequality}, Lemma \ref{L2 decay : stri-convex} and the Sobolev inequality, we arrive at
\begin{equation}
    \begin{aligned}
\begin{aligned}
\int_{\mathbb{R}} \phi_x^4(t, x) d x \leq\left\|\phi_x\right\|_{L^{\infty}}^2 {}\left\|\phi_x\right\|_{L^2}^2 & \leq C\left\|\phi_{x x}\right\|_{L^2}\left\|\phi_x\right\|_{L^2} {}\left\|\phi_x\right\|_{L^2}^2=C\left\|\phi_{x x}\right\|_{L^2}\left\|\phi_x\right\|_{L^2}^3 \\
& \leq C\left\|\phi_{x x}\right\|_{L^2}\left(\left\|\phi_{x x}\right\|_{L^2} {}\|\phi\|_{L^2}\right)^{\frac{3}{2}} \\
& \leq C\left\|\phi_{x x}\right\|_{L^2}^{\frac{5}{2}}\|\phi\|_{L^2}^{\frac{3}{2}} \leq C\left\|\phi_{x x}\right\|_{L^2}^{\frac{5}{2}}(1+t)^{\left(-\frac{1}{4}+\frac{\delta}{2}\right) \times \frac{3}{2}} \\
& \leq C(1+t)^{-\frac{3}{8}+\frac{3}{4} \delta}\left\|\phi_{x x}\right\|_{L^2}^{\frac{5}{2}},
\end{aligned}
\end{aligned}
\end{equation}
which yields
    \begin{align}
        \label{I9 esti}
I_9 & \leq \varepsilon \int_0^t(1+\tau)^\gamma\left\|\phi_{t x}(\tau)\right\|_{L^2}^2 d \tau+C \int_0^t(1+\tau)^\gamma \int_{\mathbb{R}} \phi_x^4 d x d \tau \notag\\
& \leq \varepsilon \int_0^t(1+\tau)^\gamma\left\|\phi_{t x}(\tau)\right\|_{L^2}^2 d \tau+C \int_0^t(1+\tau)^{\gamma-\frac{3}{8}+\frac{3}{4} \delta}\left\|\phi_{x x}\right\|_{L^2}^{\frac{5}{2}} d \tau \notag\\
& =\varepsilon \int_0^t(1+\tau)^\gamma\left\|\phi_{t x}(\tau)\right\|_{L^2}^2 d \tau+C \int_0^t\left\{(1+\tau)^{\frac{\gamma}{2}}\left\|\phi_{x x}\right\|_{L^2}\right\} {}\left\{(1+\tau)^{\frac{\gamma}{2}-\frac{3}{8}+\frac{3}{4} \delta}\left\|\phi_{x x}\right\|_{L^2}^{\frac{3}{2}}\right\} d \tau \notag\\
& \leq \varepsilon \int_0^t(1+\tau)^\gamma\left(\left\|\phi_{t x}\right\|_{L^2}^2+\left\|\phi_{x x}\right\|_{L^2}^2\right) d \tau+C \int_0^t(1+\tau)^{\gamma-\frac{3}{4}+\frac{3}{2} \delta}\left\|\phi_{x x}\right\|_{L^2}^3 d \tau \notag\\
& =\varepsilon \int_0^t(1+\tau)^\gamma\left(\left\|\phi_{t x}\right\|_{L^2}^2+\left\|\phi_{x x}\right\|_{L^2}^2\right) d \tau+C \int_0^t\left\{(1+\tau)^{\frac{\gamma}{2}}\left\|\phi_{x x}\right\|_{L^2}\right\} {}\left\{(1+\tau)^{\frac{\gamma}{2}-\frac{3}{4}+\frac{3}{2} \delta}\left\|\phi_{x x}\right\|_{L^2}^2\right\} d \tau \notag\\
& \leq \varepsilon \int_0^t(1+\tau)^\gamma\left(\left\|\phi_{t x}\right\|_{L^2}^2+\left\|\phi_{x x}\right\|_{L^2}^2\right) d \tau+C \int_0^t(1+\tau)^{\gamma-\frac{3}{2}+3 \delta}\left\|\phi_{x x}\right\|_{L^2}^4 d \tau \notag\\
& \leq \varepsilon \int_0^t(1+\tau)^\gamma\left(\left\|\phi_{t x}\right\|_{L^2}^2+\left\|\phi_{x x}\right\|_{L^2}^2\right) d \tau+C(1+t)^{\gamma-\frac{3}{2}+3 \delta} \int_0^t\left\|\phi_{x x}(\tau)\right\|_{L^2}^2 d \tau \notag\\
& \leq \varepsilon \int_0^t(1+\tau)^\gamma\left(\left\|\phi_{t x}\right\|_{L^2}^2+\left\|\phi_{x x}\right\|_{L^2}^2\right) d \tau+C(1+t)^{\gamma-\frac{3}{2}+3 \delta}.
\end{align}
\noindent Next, from the Sobolev inequality and Lemma \ref{Interpolation inequality}, one could deduce that 
    \begin{align}
\int_{\mathbb{R}} \phi^2 {} \phi_x^2(t, x) d x & \leq\|\phi\|_{L^{\infty}}^2 {}\left\|\phi_x\right\|_{L^2}^2 \leq C\|\phi\|_{L^2} {}\left\|\phi_x\right\|_{L^2} {}\left\|\phi_x\right\|_{L^2}^2 \notag\\
& =C\|\phi\|_{L^2} {}\left\|\phi_x\right\|_{L^2}^3 \leq C\|\phi\|_{L^2} {}\left(\left\|\phi_{x x}\right\|_{L^2} {}\|\phi\|_{L^2}\right)^{\frac{3}{2}} \notag\\
& =C\|\phi\|_{L^2}^{\frac{5}{2}}\left\|\phi_{x x}\right\|_{L^2}^{\frac{3}{2}} \leq C(1+t)^{\left(-\frac{1}{4}+\frac{\delta}{2}\right) {} \frac{5}{2}}\left\|\phi_{x x}\right\|_{L^2}^{\frac{3}{2}} \notag\\
& =C(1+t)^{-\frac{5}{8}+\frac{5}{4} \delta}\left\|\phi_{x x}\right\|_{L^2}^{\frac{3}{2}},
\end{align}
\noindent which implies that
    \begin{align}
    \label{I12 esti}
I_{12} & \leq \frac{\varepsilon}{2}\int_0^t(1+\tau)^\gamma\|\phi_{xx}(\tau)\|_{L^2}^2\,d\tau+C \int_0^t(1+\tau)^{\gamma-\frac{5}{8}+\frac{5}{4} \delta}\left\|\phi_{x x}(\tau)\right\|_{L^2}^{\frac{3}{2}} d \tau \\
& =\frac{\varepsilon}{2}\int_0^t(1+\tau)^\gamma\|\phi_{xx}(\tau)\|_{L^2}^2\,d\tau+C \int_0^t(1+\tau)^{\frac{1}{4} \gamma-\frac{5}{8}+\frac{5}{4} \delta} {}\left\{(1+\tau)^\gamma\left\|\phi_{x x}\right\|_{L^2}^2\right\}^{\frac{3}{4}} d \tau \\
& \leq \varepsilon \int_0^t(1+\tau)^\gamma\left\|\phi_{x x}(\tau)\right\|_{L^2}^2 d \tau+C \int_0^t(1+\tau)^{\gamma-\frac{5}{2}+5 \delta} d \tau \\
& \leq \varepsilon \int_0^t(1+\tau)^\gamma\left\|\phi_{x x}(\tau)\right\|_{L^2}^2 d \tau+C(1+t)^{\gamma-\frac{3}{2}+5 \delta},
\end{align}
\noindent here we have applied the Young inequality. Moreover, Lemma \ref{Properties of the smooth rarefaction wave 1}, Lemma \ref{L2 decay : stri-convex} and the mean-value formula together yield that
\begin{align}
 \label{I10,I11,I13,I14 esti}
& I_{10}+I_{11}+I_{13}+I_{14} \notag\\
\leq & \varepsilon \int_0^t(1+\tau)^\gamma\left(\left\|\phi_{t x}\right\|_{L^2}^2+\left\|\phi_{x x}\right\|_{L^2}^2\right) d \tau+C \int_0^t(1+\tau)^\gamma\left\|U_x^r\right\|_{L^{\infty}}^2\left\|\phi_x\right\|_{L^2}^2 d \tau \notag\\
& +C \int_0^t(1+\tau)^\gamma\left\|U_x^r\right\|_{L^{\infty}}^2\|\phi\|_{L^2}^2 d \tau+C \int_0^t(1+\tau)^\gamma\left(\|F(\tau)\|_{L^2}^2+\left\|F_x(\tau)\right\|_{L^2}^2\right) d \tau \notag\\
\leq & \varepsilon \int_0^t(1+\tau)^\gamma\left(\left\|\phi_{t x}\right\|_{L^2}^2+\left\|\phi_{x x}\right\|_{L^2}^2\right) d \tau+C \int_0^t(1+\tau)^{\gamma-1}\left\|\phi_x\right\|_{L^2}^2 d \tau \notag\\
& +C \int_0^t(1+\tau)^{\gamma-2}\|\phi\|_{L^2}^2 d \tau+C \int_0^t(1+\tau)^{\gamma-\frac{5}{2}} d \tau \notag\\
\leq & \varepsilon \int_0^t(1+\tau)^\gamma\left(\left\|\phi_{t x}\right\|_{L^2}^2+\left\|\phi_{x x}\right\|_{L^2}^2\right) d \tau+C(1+t)^{\gamma-\frac{3}{2}+\delta}.
\end{align}
\noindent where we have used (\ref{phi_x L^2 weighted estimate}).
Consequently, substituting (\ref{I6 esti}), (\ref{I7, I8 esti}), (\ref{I9 esti}), (\ref{I10,I11,I13,I14 esti}) and (\ref{I12 esti}) into (\ref{derivative esti-time weighted initial-stri convex}) and then taking $\varepsilon>0$ sufficiently small , it holds that that for any fixed $\delta>0$, $\gamma>5/2$, 
\begin{equation}
    \begin{aligned}
& (1+t)^\gamma\left(\left\|\phi_x(t)\right\|_{L^2}^2+\left\|\phi_{t x}(t)\right\|_{L^2}^2+\left\|\phi_{x x}(t)\right\|_{L^2}^2\right)+\int_0^t(1+\tau)^\gamma\left(\left\|\phi_{t x}(\tau)\right\|_{L^2}^2+\left\|\phi_{x x}(\tau)\right\|_{L^2}^2\right) d \tau \\
& \leq C\left(\left\|\phi_{0 x}\right\|_{L^2}^2+\left\|\phi_{0 x x}\right\|_{L^2}^2+\left\|\phi_{1 x}\right\|_{L^2}^2\right)+C(1+t)^{\gamma-\frac{3}{2}+5 \delta}.
\end{aligned}
\end{equation}
\end{proof}

Based on Lemma \ref{L2 decay : stri-convex} and \ref{derivative L2 decay: stri-convex}, it holds that for any fixed $\delta>0$, $\|\phi(t)\|_{L^2} \leq C(1+t)^{-\frac{1}{4}+\frac{1}{2}\delta},\left\|\phi_x(t)\right\|_{L^2} \leq C(1+t)^{-\frac{3}{4}+\frac{5}{2}\delta}$. Thus, owing to the Sobolev inequality,
we arrive at 
\begin{equation}
    \|\phi(t)\|_{L^{\infty}} \leq C\|\phi\|{ }_{L^2}^{1 / 2}\left\|\phi_x\right\|_{L^2}^{1 / 2} \leq C(1+t)^{-\frac{1}{2}+\frac{3}{2}\delta}.
\end{equation}
\noindent Therefore, since $\delta>0$ is arbitrary, for fixed $\eta>0$, take 
$\delta=\frac{2}{3}\eta>0$, we obtain
$\|\phi(t)\|_{L^{\infty}}\leq C(1+t)^{-\frac{1}{2}+\frac{3}{2}\delta}= C(1+t)^{-\frac{1}{2}+\eta}$. For convenience of notation, we still write $\|\phi(t)\|_{L^{\infty}} \leq C(1+t)^{-\frac{1}{2}+\delta}$, for fixed $\delta>0$.
 Moreover, according to (\ref{expression of psi_1}), we have
\begin{equation}
\label{expression to decay rate}
    \begin{aligned}
\left\|\varPsi(t)\right\|_{L^2} & \leq e^{-t}\left\|\varPsi_0\right\|_{L^2}+C \int_0^t\left(\|\phi(\tau)\|_{L^2}+\left\|\phi_x(\tau)\right\|_{L^2}+\left\|e_2(\tau)\right\|_{L^2}\right) e^{-(t-\tau)} d \tau \\
& \leq C e^{-t}+C \int_0^t(1+\tau)^{-\frac{1}{4}+\delta} e^{-(t-\tau)} d \tau \\
& \leq C e^{-t}+C(1+t)^{-\frac{1}{4}+\delta} \leq C(1+t)^{-\frac{1}{4}+\delta},
\end{aligned}
\end{equation}
\noindent and 
\begin{equation}
    \begin{aligned}
\left\|\varPsi_{ x}(t)\right\|_{L^2} & \leq C \int_0^t\left(\left\|\phi_x(\tau)\right\|_{L^2}+\left\|\phi_{x x}(\tau)\right\|_{L^2}+\left\|U_x^r(\tau)\right\|_{L^{\infty}}\|\phi(\tau)\|_{L^2}+\left\|e_{2 x}(\tau)\right\|_{L^2}\right) e^{-(t-\tau)} d \tau \\
& \qquad+e^{-t}\left\|\varPsi_{0 x}\right\|_{L^2} \\
& \leq C e^{-t}+C \int_0^t(1+\tau)^{-\frac{3}{4}+\delta} e^{-(t-\tau)} d \tau \\
& \leq C e^{-t}+C(1+t)^{-\frac{3}{4}+\delta} \leq C(1+t)^{-\frac{3}{4}+\delta},
\end{aligned}
\end{equation}
\noindent where we have used the Minkowski inequality and Lemma \ref{Properties of the smooth rarefaction wave 1}. Thus one could have \begin{equation}
    \|\psi(t)\|_{L^{\infty}} \leq C\|\psi\|{ }_{L^2}^{1 / 2}\left\|\psi_x\right\|_{L^2}^{1 / 2} \leq C(1+t)^{-\frac{1}{2}+\delta}.
\end{equation}
Finally, Lemma \ref{Properties of the smooth rarefaction wave 1}(5) gives $\|U^r(t)-u^r(\cdot/t)\|_{L^\infty}\leq C(1+t)^{-1/2}$. Since $f$ is Lipschitz continuous on the state interval, the triangle inequality yields Theorem \ref{Main Th 1}(iii). This completes the proof of Theorem \ref{Main Th 1}.

\section{The partially linearly degenerate case}
In this section, we prove Theorem \ref{Main Th 2.1}, which addresses the two-wave case corresponding to \(a<u_-<b<u_+\). As in \cite{MR4983541}, we treat the equivalent simplified case
\begin{equation}
\label{simplified assumption 2}
    \begin{cases}f^{\prime \prime}(u)>0, & u \in[0,+\infty), \\ f^{\prime \prime}(u)=0, & u \in (-\infty, 0).\end{cases}
\end{equation}
\noindent Here $f^{\prime\prime}(0)$ denotes the right derivative. We take $u_-<0<u_+$ in the following parts.

\subsection{The smooth rarefaction wave of Type 2 and the viscous contact wave}

For the partially linearly degenerate case, we first provide the construction of the smooth rarefaction wave of Type 2 and the viscous contact wave.

As for the smooth rarefaction wave II, we choose a different initial condition for the approximation than in the case of the smooth rarefaction wave I. Following the arguments in \cite{matsumura1986, matsumura2012asymptotic}, we define a different smooth approximation \(w(t, x ; w_{-}, w_{+})\) of \(w^r\left(\frac{x}{t} ; w_{-}, w_{+}\right)\) by the solution of the Cauchy problem
\begin{equation}
   \left\{\begin{array}{l}
w_t+w w_x=0,\quad t>0, x \in \mathbb{R}, \\
w(0, x)=w_0(x):=\frac{w_{-}+w_{+}}{2}+\frac{w_{+}-w_{-}}{2} \tanh x,\quad x \in \mathbb{R} .
\end{array}\right.
\end{equation}
\noindent and the same smooth approximation $U^r\left(t, x ; u_{-}, u_{+}\right)$of $u^r\left(\frac{x}{t} ; u_{-}, u_{+}\right)$by
\begin{equation*}
U^r\left(t, x ; u_{-}, u_{+}\right):=(\lambda)^{-1}\left(w\left(t, x ; \lambda_{-}, \lambda_{+}\right)\right).
\end{equation*}
 From \cite{matsumura2012asymptotic}, \cite{yoshida2014decay} and \cite{li2020}, we could obtain the following properties of $U^r(t, x)$.

\begin{lemma}[Properties of the smooth rarefaction wave of Type 2]
\label{Properties of the smooth rarefaction wave 2}
    Let $f\in C^4([u_-,u_+])$ and $f^{\prime\prime}>0$ on this interval, with one-sided derivatives at its endpoints. Then
    
(1) $U^r(t, x)$ is the $C^2$ global solution in time of the following Cauchy problem
\begin{equation*}
    \left\{\begin{array}{l}
U^r_t+f(U^r)_x=0,  \\
\left.U^r\right|_{t=0}=U^r_0(x)=\left(f^{\prime}\right)^{-1}\left(w_0(x)\right) .
\end{array}\right.
\end{equation*}

(2) $U^r_x>0$, $u_-<U^r(t, x)<u_+$.

(3)For $ 1\leq p\leq +\infty$,  there exists a constant $C_p>0$ depending on $p$ such that for $t\geq 0$,
$$
\begin{aligned}
& \left\|U_x^r(t)\right\|_{L^p},\left\|U_t^r(t)\right\|_{L^p} \leq C_p(1+t)^{-1+\frac{1}{p}}, \\
& \left\|U_{x x}^r(t)\right\|_{L^p},\left\|U_{t t}^r(t)\right\|_{L^p},\left\|U_{t x}^r(t)\right\|_{L^p},\left\|U_{x x x}^r(t)\right\|_{L^p},\left\|U_{t t x}^r(t)\right\|_{L^p} \leq C_p(1+t)^{-1}.
\end{aligned}
$$

(4) $\sup\limits _{x \in \mathbb{R}}\left|U^r(t,  x)-u^r\left(\frac{x}{t}\right)\right| \rightarrow 0$ as $t \rightarrow +\infty $.

(5){For any fixed $\delta>0$ and 
$1\leq p\leq+\infty$, there is a constant $C>0$ such that for $t>0$
\begin{equation*}
  \left\|U^r(t, x)-u^r\left(\frac{x}{t}\right)\right\|_{L^p{}} \leqslant C_{p, \delta}(1+t)^{(-1+\delta)\left(1-\frac{1}{p}\right)} .
\end{equation*}}
\noindent Moreover, for $t\geq 0$,
\begin{equation*}
    \left\|U^r(t, x)-u^r\left(\frac{x}{1+t}\right)\right\|_{L^p{}} \leqslant C_{p, \delta}(1+t)^{(-1+\delta)\left(1-\frac{1}{p}\right)}.
\end{equation*}
\end{lemma}
\begin{remark}
In Sections 3 and 4, we employ two distinct types of smooth rarefaction waves. The rationale is as follows. The Type 1 smooth rarefaction wave possesses a sharper decay rate for high-order derivatives, which enables us to establish the optimal decay rate of the perturbation in the strictly convex flux case. However, the \(L^\infty\) decay rate of \(U^r(t,x) - u^r(\frac{x}{t})\) for Type 1 is insufficient to derive the interaction estimates of nonlinear waves, such as those in Lemma \ref{X(t) properties} and Lemma \ref{F(U,U^r) properties}. In view of (5) in Lemma \ref{Properties of the smooth rarefaction wave 2}, we adopt the Type 2 approximation for the partially linearly degenerate flux.
    \end{remark}

Next, we will concentrate on the properties of the viscous contact wave $U(t,x)$. 
\begin{lemma}[Properties of the viscous contact wave \cite{MR4983541}]
    \label{Properties of the viscous contact wave}
Suppose that the far field states satisfy $u_-<u_+$, then

(1) Defined by (\ref{Expression of VC}) with $\mu=1-\bar\lambda^2>0$, $U(t,x)\in C^{\infty} ((0, +\infty) \times \mathbb{R}) $ is one self-similar solution of the below Cauchy problem 
\begin{equation*}
    \begin{aligned}
\begin{cases} U_t+\bar{\lambda}U_x-(1-\bar{\lambda}^2)  U_{xx}=0,\qquad t>0, \quad x \in \mathbb{R}, \\ U(0, x)= \begin{cases}u_{-}, & x<0, \\ u_{+}, & x>0,\end{cases} \\ \lim\limits _{x \rightarrow \pm \infty} U(t, x)=u_{ \pm},\qquad t \geq 0.\end{cases}
\end{aligned}
\end{equation*}
 
(2) $u_{-}<U(t, x)<u_{+} $ and $  U_x(t, x)>0$ for $t>0$ $x \in \mathbb{R}$.

(3) If $\bar{\lambda}=0$, for $1 \leq p \leq \infty $, $k,l\in\mathbb{N}_0$, $k+l\geq1$, there exists a  positive constant $C_{p,k,l} $, such that for  $t>0$,

\begin{equation}
\left\|\partial_x^k \partial_t^l U(t)\right\|_{L^p} \leq C_{p, k, l} t^{-\frac{k}{2}-l+\frac{1}{2 p}}.
\end{equation}

(4) If $\bar{\lambda}\neq0$,
for $1 \leq p \leq \infty $, $k,l\in\mathbb{N}_0$, $k+l\geq1$, there exists a  positive constant $ C_{p,k,l}$, such that for $t>0$,
\begin{equation}
\left\|\partial_x^k \partial_t^l U(t)\right\|_{L^p} \leq C_{p, k, l} t^{-\frac{k+l}{2}+\frac{1}{2 p}}.
\end{equation}

(5) $\lim\limits _{t \rightarrow +\infty} \sup\limits _{x \in \mathbb{R}}|U(1+t, x)-U(t, x)|=0$.
\end{lemma}

Based on (\ref{simplified assumption 2}), Lemma \ref{Properties of the smooth rarefaction wave 2} and Lemma \ref{Properties of the viscous contact wave}, we present the construct of the composite wave:
\begin{equation*}
    U_{\mathrm {multi }}(t, x):=U\left(t,x ; u_{-}, 0\right)+u^{r}\left(\frac{x}{t} ; 0, u_{+}\right).
\end{equation*}
\noindent And define the smooth asymptotic state
\begin{equation*}
    \widetilde{U}(t,x):=U^r(t,x;0,u_+)+U(1+t,x;u_-,0).
\end{equation*}
\noindent Now take the perturbation 
$$\phi(t, x):=u(t, x)-\widetilde{U}(t, x), \psi(t, x):=v(t, x)-f({\widetilde{U}}(t, x)),$$
\noindent with the initial data $(\phi_0, \psi_0)(x)$. The uniform energy estimate for the perturbation have been obtained in \cite{MR4983541}.
\begin{lemma}[Uniform energy estimate \textrm{II}\cite{MR4983541}]
\label{uniform esti: degenerate case}
    It holds that there exists a constant $C\left(\phi_0, \psi_0\right)>0$ depending on $\left\|\phi_0\right\|_{H^2},\left\|\psi_0\right\|_{H^2}$ and on the fixed flux and wave profiles, such that for each $0<t<+\infty$, the unique global smooth solution $\left(u(t, x), v(t, x)\right)$ in 
Lemma \ref{Global existence} satisfies
\begin{equation*} \|\phi(t)\|_{H^2}^2+\left\|\phi_t(t)\right\|_{H^1}^2+\int_0^t G(\tau)+\left\|\phi_x(\tau)\right\|_{H^1}^2+\left\|\phi_t(\tau)\right\|_{H^1}^2\mathrm{d} \tau \leq C\left(\phi_0, \psi_0\right) ,
\end{equation*}
\noindent and
\begin{equation*} \|\psi(t)\|_{H^2}^2+\left\|\psi_t(t)\right\|_{H^1}^2+\int_0^t\left\|\psi_x(\tau)\right\|_{H^1}^2+\left\|\psi_t(\tau)\right\|_{H^1}^2\mathrm{d} \tau \leq C\left(\phi_0, \psi_0\right),
\end{equation*}
\noindent where
\begin{equation*}
    \begin{aligned}
G(t):= & \left(\int_{\widetilde{U}+\phi \geq 0, \widetilde{U} \geq 0} \phi^2 \widetilde{U}_x \mathrm{~d} x\right)(t)+\left(\int_{\widetilde{U}+\phi \geq 0, \widetilde{U}<0}(\widetilde{U}+\phi)^2 \widetilde{U}_x \mathrm{~d} x\right) \\
& +\left(\int_{\widetilde{U}+\phi<0, \widetilde{U} \geq 0} \widetilde{U}^2 \widetilde{U}_x \mathrm{~d} x\right)(t) .
\end{aligned}
\end{equation*}

\end{lemma}

The corresponding main theorem is as follows.
\begin{theorem}[Partially linearly degenerate case]
    \label{Main Th 2.1}
Suppose that flux $f(u) \in C^1(\mathbb{R}) \cap C^4([0, +\infty))$ satisfy (\ref{simplified assumption 2}) and the sub-characteristic condition (\ref{sub-characteristic}), the far field states satisfy $u_-<0<u_+$. If $u_0-u_0^R, v_0-f\left(u_0^R\right) \in L^1(\mathbb{R}) \cap L^2(\mathbb{R}) $ and the weak derivative $u_{0 x}, v_{0 x} \in H^1(\mathbb{R})$, then the Cauchy problem
(\ref{Original equation: Jin-Xin system})-(\ref{initial values}) has a unique global smooth solution $(u(t,x), v(t,x))$, satisfying
\begin{equation*}
\sup_{t\geq0}\|u(t)\|_{L^\infty}\leq2N_0,\qquad \sup_{t\geq0}\|v(t)\|_{L^\infty}\leq B(N_0),
\end{equation*}
\noindent where the definition of $B(N_0)$ is the same as in Theorem \ref{Main Th 1}. And 
for any fixed $\delta>0$, there exists a constant $C>0$, 
such that for $t>0$,
\begin{equation*}
    \left\|\left(u(t, x)-U_{\mathrm{multi}}(t, x), v(t, x)-f\left(U_{\mathrm{multi}}(t, x)\right)\right)\right\|_{L^1{}} \leq C (1+t)^\delta,
\end{equation*}
\noindent Moreover, it holds that for the above \(\delta>0\) , there exists  constants $T>0$ and $C>0$, for all $t\geq T$, such that
\begin{equation*}
\label{decay rate: degenerate case}
 \left\|\left(u(t, x)-U_{\mathrm{multi}}(t, x), v(t, x)-f\left(U_{\mathrm{multi}}(t, x)\right)\right)\right\|_{L^{\infty}{}} \leq C (1+t)^{-\frac{1}{3}+\delta},
\end{equation*}
\noindent where, 
$$
U_{\mathrm {multi }}(t, x):=U\left(t,x ; u_{-}, 0\right)+u^{r}\left(\frac{x}{t} ; 0, u_{+}\right).
$$
\end{theorem}

\subsection{The interaction of the nonlinear waves}
 In the interaction and perturbation estimates below, an unsubscripted argument $U$ or a derivative of $U$ denotes the contact profile evaluated at $(1+t,x)$, unless its original time argument is displayed explicitly. We first transform the nonlinear flux. Let $g(u)=f(u)-\bar{\lambda} u-f(0)$, which leads to $g^{\prime}(u)=f^{\prime}(u)-\bar{\lambda}$. 
The basic properties of $g(u)$ are as follows:
\begin{lemma}
\label{properties of $g$}
$g(u)=f(u)-\bar{\lambda} u-f(0)$ has the following properties:

(1) $g(u) \in C^1(\mathbb{R}) \cap C^4([0,+\infty))$;

(2) For $u > 0$, $g^{\prime \prime}(u)>0$;

(3) For $u < 0$, $g(u)=0$;

(4) $g(0)=g^{\prime}(0)=0$, and the second-order right derivative of $g(u)$ at $u=0$
is positive;

(5) For a fixed finite number $A>0$, there exists a constant $C>0$, such that for  $0\leq u\leq A$,
\begin{equation*}
    0\leq g(u) \leq C u^2, 0\leq g^{\prime}(u) \leq C u.
\end{equation*}

(6)\(g'\) is Lipschitz continuous on any bounded interval \(I \subset \mathbb{R}\). In particular, there exists a constant \(C > 0\) such that
\begin{equation*}
    \left|g^{\prime}(\widetilde{U}+\phi)-g^{\prime}(\widetilde{U})\right| \leq C|\phi|.
\end{equation*}
\end{lemma}
\begin{proof}
It suffices to establish (6); for the proofs of the remaining items, we refer to \cite{MR4983541}.  We now prove (6). Indeed, for $u_1,u_2\in\mathbb{R}$, if $u_1,u_2\geq0$ or $u_1,u_2<0$, the conclusion follows immediately. If $u_1\geq0$ and $u_2<0$, it holds from (5) that $|g'(u_1)-g'(u_2)|=g'(u_1)\leq Cu_1\leq C(u_1-u_2)$. 
The remaining case follow by symmetry.
\end{proof}
 Next, according to Lemma \ref{Properties of the smooth rarefaction wave 2} and Lemma \ref{Properties of the viscous contact wave}, 
we can rewrite the problem \eqref{damped wave equ} in terms of $\phi$ as
\begin{equation}
\label{perturbation equ: degenerate case}
\left\{\begin{array}{l}
\phi_{t t}+\phi_t+(g(\widetilde{U}+\phi)-g(\widetilde{U}))_x+\bar{\lambda} \phi_x-\phi_{x x}=-U_{t t}^r-U_{t t}+\bar{\lambda}^2 U_{x x}+F\left(U, U^r\right), \\
\left.\phi\right|_{t=0}=u_0(x)-\widetilde{U}(0, x) \triangleq \phi_0(x), \\
\left.\phi_t\right|_{t=0}=u_1(x)-\widetilde{U}_t(0, x) \triangleq \phi_1(x) .
\end{array}\right.
\end{equation}
\noindent Here
$$
F(U,U^r):=-\left(g^{\prime}\left(U+U^r\right)-g^{\prime}\left(U^r\right)\right)  U^r_x-g^{\prime}\left(U+U^r\right)  U_x+U_{x x}^r.
$$

As the preparations for the estimate of the time-weighted $L^2$ estimate, it is necessary to give the basic properties of $F(U,U^r)$.
Based on the properties of $g(u)$, to get the $L^p$ estimate of $F(U,U^r)$, it is  natural to divide the integral region $\mathbb{R}$ depending on the sign of 
$\widetilde{U}=U^r +U$. 

From $\widetilde{U}_x(t,x)=U^r_x(t,x)+U_x(1+t,x)>0$ and $\lim\limits_{x \rightarrow{-\infty}}\widetilde{U}(t,x)=u_-<0$, $\lim\limits_{x \rightarrow{+\infty}}\widetilde{U}(t,x)=u_+>0$, we can have that for each fixed $t\geq0$, there exists an unique $x=X(t)$ such that 
\begin{equation}
\label{X(t)}
    \widetilde{U}(t,X(t))={U}^r(t,X(t))+{U}(1+t,X(t))=0,
\end{equation}
\noindent which is equivalent to:
\begin{equation}
\label{relation: Ur and U}
    U^r(t, X(t))=-U(1+t, X(t))=\frac{\left|u_{-}\right|}{\sqrt{\pi}} \int_{\frac{X(t)-\bar{\lambda}(1+t)}{\sqrt{4 \mu(1+t)}}}^{+\infty} \mathrm{e}^{-\xi^2} \mathrm{~d} \xi>0 .
\end{equation}
And, the estimate of $X(t)$ is as follow.
\begin{lemma}[\cite{MR4983541}]
\label{X(t) properties}
$X(t)$ defined by (\ref{X(t)}) has the following properties.

(1) There exists constants $C>0$ and $T_0>0$ such that for any $t\geq T_0$, such that
    $$
\bar{\lambda}(1+t)+\sqrt{4(1-\bar{\lambda}^2)(1+t)}\leq X(t)\leq \lambda_+(1+t).
    $$

(2) There exist constants $C>0$ and $T_0>0$ such that for  $t\geq T_0$ and $\epsilon\in (0,1)$,
    $$
\begin{aligned}
& 0<\frac{X(t)-\bar{\lambda}(1+t)}{1+t} \leq C(1+t)^{-1+\epsilon}+C \int_{\frac{X(t)-\bar{\lambda}(1+t)}{\sqrt{4\left(1-\bar{\lambda}^2\right)(1+t)}}}^{+\infty} \mathrm{e}^{-\xi^2} \mathrm{~d} \xi, \\
& \frac{X(t)-\bar{\lambda}(1+t)}{1+t} \geq C \int_{\frac{X(t)-\bar{\lambda}(1+t)}{\sqrt{4\left(1-\bar{\lambda}^2\right)(1+t)}}}^{+\infty} \mathrm{e}^{-\xi^2} \mathrm{~d} \xi-C(1+t)^{-1+\epsilon}.
\end{aligned}
    $$

(3) Denote that $Y(t)=\frac{X(t)-\bar{\lambda}(1+t)}{\sqrt{4\left(1-\bar{\lambda}^2\right)(1+t)}}$, then there exists a constant $C>0$ such that for $t\geq T_0$
$$
Y(t) \leq\left(C+\ln (1+t)^{\frac{1}{2}}\right)^{\frac{1}{2}} .
$$

(4) For any $0<\delta<1$, there exist constants  $C_\delta>0$ and $T_0>0$, such that for any $t\geq T_0$
    $$
    \left(C_\delta+\ln (1+t)^{\frac{1}{2(1+\delta)}}\right)^{\frac{1}{2}}  \leq Y(t).
    $$

(5) There exists a constant $C>0$ such that for $t\geq 0$
    $$
    0<U^r(t, X(t))  \leq C \sqrt{\frac{\ln (2+t)}{1+t}}.
    $$

\end{lemma}
According to Lemma \ref{X(t) properties}, we have the following preliminary estimates for the residual \(F(U, U^r)\), which can be found in \cite{MR4983541}.

\begin{lemma}[Lemmas 3.5--3.6 in \cite{MR4983541}]
    \label{F(U,U^r) properties}
For any $(\delta, p) \in(0, \infty) \times[1, \infty)$, there exists a constant $C$  such that for $t\geq 0$, the following holds that
$$
\left\|F\left(U, U^r\right)(t)\right\|_{L^p} \leq C(1+t)^{-\frac{1}{2}\left(1+\frac{1}{p}\right)+\frac{\delta}{2}},\qquad \left\|\partial_x F\left(U, U^r\right)(t)\right\|_{L^p} \leq C(1+t)^{-1+\frac{\delta}{2}}.
$$
\end{lemma}

\subsection{The $L^1$ estimate of $(\phi, \psi)$}

Next, we still derive the $L^1$ estimate by the similar way in the previous section.
\begin{lemma}
\label{L1 esti: degenerate case}
Under the same assumptions as in Theorem \ref{Main Th 2.1}, for $\delta>0$  as in Lemma \ref{F(U,U^r) properties}, 
there exists a constant $C_\delta>0$, 
such that for $t\geq0$,
\begin{equation*}
\left\|\left(\phi(t, x), \psi(t, x)\right)\right\|_{L^1{}} \leq C_\delta (1+t)^\frac{\delta}{2}.
\end{equation*}

\end{lemma}
\begin{proof}
We still use the similar decomposition argument to prove this Lemma. Notice that the viscous contact wave $U(t,x)$ satisfies $U_t+\bar{\lambda}U_x=(1-\bar{\lambda}^2)U_{xx}$, take $V:=f(U)-\left(1-\bar{\lambda}^2\right) U_x-\bar{\lambda}\left(1-\bar{\lambda}^2\right) U_{x x}+\left(1-\bar{\lambda}^2\right) U_{t x}$, one could get
\begin{equation}
    \begin{aligned}
U_t+V_x & =U_t+f(U)_x-\left(1-\bar{\lambda}^2\right) U_{x x}-\bar{\lambda}\left(1-\bar{\lambda}^2\right) U_{x x x}+\left(1-\bar{\lambda}^2\right) U_{t x x} \\
& =\left(U_t+\bar{\lambda} U_x-\left(1-\bar{\lambda}^2\right) U_{x x}\right)-\bar{\lambda}\left(1-\bar{\lambda}^2\right) U_{x x x}+\left(1-\bar{\lambda}^2\right) U_{t x x} \\
& =-\bar{\lambda}\left(1-\bar{\lambda}^2\right) U_{x x x}+\left(1-\bar{\lambda}^2\right) U_{t x x}=: e_3, 
\end{aligned}
\end{equation}
\noindent where we have used the facts that $U<0$ and $f^\prime(u)=f^\prime(0)=\bar{\lambda}$ on $(-\infty, 0]$. Moreover,  
\begin{equation}
    \begin{aligned}
& V_t+U_x-(f(U)-V) \\
= & \bar{\lambda} U_t-\left(1-\bar{\lambda}^2\right) U_{t x}-\bar{\lambda}\left(1-\bar{\lambda}^2\right) U_{x x t}+\left(1-\bar{\lambda}^2\right) U_{x t t} \\
& +\bar\lambda^2 U_x-\bar\lambda(1-\bar\lambda^2)U_{xx}+(1-\bar\lambda^2)U_{tx} \\
= & \bar{\lambda}\left(-\bar{\lambda} U_x+\left(1-\bar{\lambda}^2\right) U_{x x}\right)-\left(1-\bar{\lambda}^2\right) U_{t x}-\bar{\lambda}\left(1-\bar{\lambda}^2\right) U_{x x t}+\left(1-\bar{\lambda}^2\right) U_{x t t} \\
& +U_x-f(U)+f(U)-\left(1-\bar{\lambda}^2\right) U_x-\bar{\lambda}\left(1-\bar{\lambda}^2\right) U_{x x}+\left(1-\bar{\lambda}^2\right) U_{t x} \\
= & -\bar{\lambda}\left(1-\bar{\lambda}^2\right) U_{x x t}+\left(1-\bar{\lambda}^2\right) U_{x t t}=: e_4.
\end{aligned}
\end{equation}
\noindent That is the approximate Jin-Xin model with respect to $(U,V)$
\begin{equation}
\label{Jin-Xin L^1 eq2.1}
    \left\{\begin{aligned}
U_t+V_x & =e_3, \\
V_t+U_x & =f(U)-V+e_4.
\end{aligned}\right.
\end{equation}
 \noindent Take $\widetilde{V}:=V^r+V=f\left(U^r\right)+f(U)-\left(1-\left(f^{\prime}\left(U^r\right)\right)^2\right) U_x^r-\left(1-\bar{\lambda}^2\right) U_x-\bar{\lambda}\left(1-\bar{\lambda}^2\right) U_{x x}+\left(1-\bar{\lambda}^2\right) U_{t x}$, thus by adding  (\ref{Jin-Xin L^1 eq1.3})
 and (\ref{Jin-Xin L^1 eq2.1}), we could obtain
\begin{equation}
\label{middle state of L^1}
\begin{cases}
\widetilde {U}_t + \widetilde {V}_x = e_5, \\
\widetilde {V}_t + \widetilde {U}_x = \left( f(U^r) + f(U) \right) - \widetilde {V} + e_6.
\end{cases} 
\end{equation}
\noindent Here $e_5:=e_1+e_3$, $e_6:=e_2+e_4$. 
In view of Lemma \ref{Properties of the smooth rarefaction wave 2} and Lemma \ref{Properties of the viscous contact wave}, one could have $\left\|e_5(t)\right\|_{L^1},\left\|e_6(t)\right\|_{L^1} \leq C(1+t)^{-1}$. To get the the approximate Jin-Xin model with respect to $\widetilde U$, we should decompose further. 
Let the correction term 
\begin{equation}
\label{correction term}
    \begin{aligned}
& Q:=f(\widetilde{U})-\left(f\left(U^r\right)+f(U)\right), \\
& \mathcal{\widetilde{V}}:=\widetilde{V}+Q=f(\widetilde{U})-\left(1-\left(f^{\prime}\left(U^r\right)\right)^2\right) U_x^r-\left(1-\bar{\lambda}^2\right) U_x-\bar{\lambda}\left(1-\bar{\lambda}^2\right) U_{x x}+\left(1-\bar{\lambda}^2\right) U_{t x}.
\end{aligned}
\end{equation}
\noindent Hence, substituting (\ref{correction term}) into (\ref{middle state of L^1}), one could have
\begin{equation}
\begin{aligned}
& \widetilde{U}_t+\widetilde{\mathcal{V}}_x=\widetilde{U}_t+\widetilde{V}_x+Q_x=e_5+Q_x=: e_7, \\
& \widetilde{\mathcal{V}}_t+\widetilde{U}_x-(f(\widetilde{U})-\widetilde{\mathcal{V}}) \\
= & \widetilde{V}_t+Q_t+\widetilde{U}_x-f(\widetilde{U})+\widetilde{V}+Q \\
= & \left(f\left(U^r\right)+f(U)-\widetilde{V}+e_6\right)+Q_t-f(\widetilde{U})+\widetilde{V}+f(\widetilde{U})-\left(f\left(U^r\right)+f(U)\right)=e_6+Q_t=: e_8.
\end{aligned}
\end{equation}
\noindent Thus the approximate Jin-Xin model with respect to $(\widetilde U, \mathcal{\widetilde{V}})$ is
\begin{equation}
\label{Jin-Xin L^1 eq2.2}
  \left\{\begin{array}{l}
\widetilde{U}_t+\mathcal{\widetilde{V}}_{ x}=e_7, \\
\mathcal{\widetilde{V}}_{1 t}+\widetilde{U}_x=f(\widetilde{U})-\mathcal{\widetilde{V}}+e_8 .
\end{array}\right.
\end{equation}

\noindent Let $\widetilde\varPsi:=v-\mathcal{\widetilde{V}}=\psi+\left(1-\left(f^{\prime}\left(U^r\right)\right)^2\right) U_x^r+\left(1-\bar{\lambda}^2\right) U_x+\bar{\lambda}\left(1-\bar{\lambda}^2\right) U_{x x}-\left(1-\bar{\lambda}^2\right) U_{t x}$; thus, via $U^r_x, U_x, U_{xx}, U_{tx}\in L^\infty(L^1(\mathbb{R}))$, $\psi_{0}(x):=\psi|_{t=0}\in L^1(\mathbb{R})$ leads to $\widetilde\varPsi_0(x):=\widetilde\varPsi|_{t=0}\in L^1(\mathbb{R})$. It holds  from (\ref{Original equation: Jin-Xin system}) and (\ref{Jin-Xin L^1 eq2.2}) that 
\begin{equation}
\label{Jin-Xin L^1 eq2.3}
    \left\{\begin{array}{l}
\phi_t+\widetilde\varPsi_{ x}=-e_7, \\
\widetilde\varPsi_{ t}+\phi_x=f(\widetilde{U}+\phi)-f(\widetilde{U})-\widetilde\varPsi-e_8 .
\end{array}\right.
\end{equation}

\noindent By applying the same $L^1$- contraction estimate argument as in the proof of Lemma \ref{L1 estimate: Stri-convex}, there remains to verify that for $\delta>0$ as in Lemma \ref{F(U,U^r) properties}, $\left\|e_7(t)\right\|_{L^1},\left\|e_8(t)\right\|_{L^1} \leq C(1+t)^{-1+\frac{\delta}{2}}$. Thus we only need to show
$\left\|Q_t(t)\right\|_{L^1},\left\|Q_x(t)\right\|_{L^1} \leq C(1+t)^{-1+\frac{\delta}{2}}$. Note that
\begin{align}
Q_x & =f^{\prime}\left(U^r+U\right)\left(U_x^r+U_x\right)-f^{\prime}\left(U^r\right) U_x^r-f^{\prime}(U) U_x \notag\\
& =\left(f^{\prime}\left(U^r+U\right)-f^{\prime}\left(U^r\right)\right) U_x^r+\left(f^{\prime}\left(U^r+U\right)-f^{\prime}(U)\right) U_x \notag\\
& =\left(g^{\prime}\left(U+U^r\right)-g^{\prime}\left(U^r\right)\right) U_x^r+g^{\prime}\left(U+U^r\right) U_x \notag\\
& =-F+U_{x x}^r, \notag\\
Q_t & =\left(f^{\prime}\left(U^r+U\right)-f^{\prime}\left(U^r\right)\right) U_t^r+\left(f^{\prime}\left(U^r+U\right)-f^{\prime}(U)\right) U_t \notag\\
& =\left(f^{\prime}\left(U^r+U\right)-f^{\prime}\left(U^r\right)\right)\left(-f^{\prime}\left(U^r\right) U_x^r\right)+\left(f^{\prime}\left(U^r+U\right)-f^{\prime}(U)\right)\left(-\bar{\lambda} U_x+\left(1-\bar{\lambda}^2\right) U_{x x}\right) \notag\\
& =-\bar\lambda Q_x-g'(U^r)\bigl(g'(\widetilde U)-g'(U^r)\bigr)U_x^r
+(1-\bar\lambda^2)g'(\widetilde U)U_{xx}.
\end{align}
\noindent In view of the proof of Lemma 3.6 in \cite{MR4983541}, since  
as $x>X(t)>\bar{\lambda}(1+t)$,
$$
\begin{gathered}
U_x(1+t, x)=\frac{\left|u_{-}\right|}{\sqrt{\pi}} {} \frac{1}{\sqrt{4\left(1-\bar{\lambda}^2\right)(1+t)}} \exp \left(-\frac{(x-\bar{\lambda}(1+t))^2}{4\left(1-\bar{\lambda}^2\right)(1+t)}\right), \\
U_{x x}(1+t, x)=-\frac{\left|u_{-}\right|}{\sqrt{\pi}} {} \frac{x-\bar{\lambda}(1+t)}{2\left(1-\bar{\lambda}^2\right)(1+t) \sqrt{4\left(1-\bar{\lambda}^2\right)(1+t)}} \exp \left(-\frac{(x-\bar{\lambda}(1+t))^2}{4\left(1-\bar{\lambda}^2\right)(1+t)}\right)<0.
\end{gathered}
$$
\noindent  So for fixed $t>0$, $U_x(1+t,x)\rightarrow{0}$, as $x\rightarrow{\pm \infty} $. Thus,  the integration by parts implies that for the above $\delta>0$,
\begin{align}
    \label{one residual esti}
\|g^{\prime}(\widetilde U)U_{xx}(1+t,x)\|_{L^1} & = \int_{\mathbb{R}} g^{\prime}\left(U+U^r\right)\left|U_{x x}\right| d x=\int_{X(t)}^{+\infty} g^{\prime}\left(U+U^r\right)\left(-U_{x x}\right) d x \notag\\
& =-\int_{X(t)}^{+\infty} g^{\prime}\left(U+U^r\right) d U_x \notag\\
& =\int_{X(t)}^{+\infty} U_x g^{\prime \prime}\left(U+U^r\right)\left(U_x+U_x^r\right) d x \notag\\
& \leq\left(\left\|U_x^r\right\|_{L^{\infty}}+\left\|U_x\right\|_{L^{\infty}}\right) \int_{X(t)}^{+\infty} g^{\prime \prime}\left(U+U^r\right) U_x d x \notag\\
& \leq C(1+t)^{-\frac{1}{2}} \int_{X(t)}^{+\infty} U_x d x=-C(1+t)^{-\frac{1}{2}} U(1+t, X(t)) \notag\\
& \leq C(1+t)^{-1+\frac{\delta}{2}},
\end{align}
\noindent where we have used (5) in Lemma \ref{X(t) properties}. This calculation is for $t\geq T_0$; the bounded interval $0\leq t\leq T_0$ is covered by smoothness of the time-shifted profiles.
To estimate the remaining term in $Q_t$, use the local Lipschitz continuity of $g'$ and split at $X(t)$:
\begin{equation*}
\begin{aligned}
E(t)&:=\int_{\mathbb R}g'(U^r)\left|g'(\widetilde U)-g'(U^r)\right|U_x^r\,dx\\
&\leq C\bigl(U^r(t,X(t))\bigr)^3
+C\|U_x^r(t)\|_{L^\infty}\int_{X(t)}^\infty |U(1+t,x)|\,dx.
\end{aligned}
\end{equation*}
Indeed, on $x<X(t)$ one has $g'(\widetilde U)=0$ and $g'(U^r)\leq C U^r$; on $x>X(t)$ use $|g'(\widetilde U)-g'(U^r)|\leq C|U|$.
The Gaussian formula and (4) in Lemma \ref{X(t) properties} give, for each $\delta>0$ and sufficiently large $t$,
\begin{equation*}
\int_{X(t)}^\infty |U(1+t,x)|\,dx
\leq C\sqrt{1+t}\,e^{-Y(t)^2}\leq C_\delta(1+t)^\frac{\delta}{2}.
\end{equation*}
Together with (5) of that lemma, this yields $E(t)\leq C_\delta(1+t)^{-1+\frac{\delta}{2}}$, after enlarging the constant on a bounded time interval. Thus the exact formula for $Q_t$, (\ref{one residual esti}), Lemma \ref{Properties of the smooth rarefaction wave 2}, and Lemma \ref{F(U,U^r) properties} imply
\begin{equation}
\left\|Q_t(t)\right\|_{L^1},\left\|Q_x(t)\right\|_{L^1} \leq C(1+t)^{-1+\frac{\delta}{2}}.
\end{equation}
The remaining proof is the same as in the proof of Lemma \ref{L1 estimate: Stri-convex}. The \(L^1\) bound in Lemma \ref{L1 esti: degenerate case} differs from that in Lemma \ref{L1 estimate: Stri-convex}, due to the different \(L^1\) decay rates of their residuals. Specifically, the residual decays as \((1+t)^{-1+\frac{\delta}{2}}\) in the former case, but as \((1+t)^{-1}\) in the latter. Therefore, we conclude the proof. 
\end{proof}

\subsection{$L^2$ decay for $(\phi, \psi)$}
\begin{lemma}
    \label{L2 decay: degenerate case}
Under the same assumptions as in Theorem \ref{Main Th 2.1}, for $\delta>0$  as in Lemma \ref{F(U,U^r) properties} and $\gamma>2$,
there exists a constant $C>0$, 
such that for $t\geq0$,
\begin{equation*}
 \begin{aligned}
& (1+t)^\gamma\left(\|\phi(t)\| _{L^2}^2+\left\|\phi_t(t)\right\|_{L^2}^2+\left\|\phi_x(t)\right\|_{L^2}^2\right)+\int_0^t(1+\tau)^\gamma\left(\left\|\phi_t(\tau)\right\|_{L^2}^2+\left\|\phi_x(\tau)\right\|_{L^2}^2+G(\tau)\right) d \tau \\
& \leq C\left(\left\|\phi_0\right\|_{L^2}^2+\left\|\phi_{0 x}\right\|_{L^2}^2+\left\|\phi_1\right\|_{L^2}^2\right)+C(1+t)^{\gamma-\frac{1}{2}+\frac{5}{4}\delta},
\end{aligned}
\end{equation*}
\noindent where
\begin{equation*}
    \begin{aligned}
G(t):= & \left(\int_{\widetilde{U}+\phi \geq 0, \widetilde{U} \geq 0} \phi^2 \widetilde{U}_x \mathrm{~d} x\right)(t)+\left(\int_{\widetilde{U}+\phi \geq 0, \widetilde{U}<0}(\widetilde{U}+\phi)^2 \widetilde{U}_x \mathrm{~d} x\right) \\
& +\left(\int_{\widetilde{U}+\phi<0, \widetilde{U} \geq 0} \widetilde{U}^2 \widetilde{U}_x \mathrm{~d} x\right)(t) .
\end{aligned}
\end{equation*}
\end{lemma}

\begin{proof}
The proof of Lemma \ref{L2 decay: degenerate case} is essentially due to He, Huang and Wang\cite{MR4983541}. Multiplying (\ref{perturbation equ: degenerate case}) by $\phi+2\phi_t$,
\begin{equation}
\label{div form 1: degenerate case}
   \begin{aligned}
& \left(\frac{1}{2} \phi^2+\phi \phi_t+\phi_t^2+\phi_x^2\right)_t+\left(\phi_t^2+\phi_x^2+2\left(\bar{\lambda}+g^{\prime}(u)\right) \phi_t \phi_x\right)+\left(\frac{\bar{\lambda}}{2} \phi^2-\phi \phi_x-2 \phi_t \phi_x\right)_x \\
& \quad+\phi {}(g(\widetilde{U}+\phi)-g(\widetilde{U}))_x+2 \phi_t {}\left(g^{\prime}(\widetilde{U}+\phi)-g^{\prime}(\widetilde{U})\right) \widetilde{U}_x \\
& =-\phi {} U_{t t}^r-\phi {}\left(U_{t t}-\bar{\lambda}^2 U_{x x}\right)-2 \phi_t {} U_{t t}^r-2 \phi_t {} U_{t t}+2 \bar{\lambda}^2 \phi_t {} U_{x x}+\left(\phi+2 \phi_t\right) {} F,
\end{aligned}
\end{equation}
\noindent where we used the fact that  $2 \phi_t {}(g(\widetilde{U}+\phi)-g(\widetilde{U}))_x=2 \phi_t {}\left(g^{\prime}(\widetilde{U}+\phi)-g^{\prime}(\widetilde{U})\right) \widetilde{U}_x+2 \phi_t {} g^{\prime}(\widetilde{U}+\phi) {} \phi_x$.
\noindent By utilizing  the same argument as in (\ref{positive-quard: stri-convex}), we arrive at
\begin{equation}
    \label{positive quara: degenerate case}
\phi_t^2+\phi_x^2+2\left(\bar{\lambda}+g^{\prime}(u)\right) \phi_t \phi_x \geq \phi_t^2+\phi_x^2-2 k\left|\phi_t\right|\left|\phi_x\right| \geq \alpha\left(\phi_t^2+\phi_x^2\right),
\end{equation}
\noindent here $\alpha$ is a positive constant.
\noindent Moreover, a direct computation 
yields that
\begin{equation}
\label{calculation of disspation: degenerate case}
    \begin{aligned}
& \phi {}(g(\widetilde{U}+\phi)-g(\widetilde{U}))_x \\
= & \{\phi {}(g(\widetilde{U}+\phi)-g(\widetilde{U}))\}_x-\left\{\int_0^\phi g(\widetilde{U}+\eta)-g(\widetilde{U}) d \eta\right\}_x+\widetilde{U}_x\left(g(\widetilde{U}+\phi)-g(\widetilde{U})-\phi {} g^{\prime}(\widetilde{U})\right) .
\end{aligned}
\end{equation}
\noindent Plugging (\ref{positive quara: degenerate case}) and \ref{calculation of disspation: degenerate case} into (\ref{div form 1: degenerate case}), we deduce that 
    \begin{align}
    \label{div form: degenerate}
& \left(\frac{1}{2} \phi^2+\phi \phi_t+\phi_t^2+\phi_x^2\right)_t+\left(\frac{\bar{\lambda}}{2} \phi^2-\phi \phi_x-2 \phi_t \phi_x+\phi {}(g(\widetilde{U}+\phi)-g(\widetilde{U}))-\int_0^\phi g(\widetilde{U}+\eta)-g(\widetilde{U}) d \eta\right)_x \notag\\
& \quad+\alpha\left(\phi_t^2+\phi_x^2\right)+\widetilde{U}_x\left(g(\widetilde{U}+\phi)-g(\widetilde{U})-\phi {} g^{\prime}(\widetilde{U})\right)+2 \phi_t {}\left(g^{\prime}(\widetilde{U}+\phi)-g^{\prime}(\widetilde{U})\right) \widetilde{U}_x \notag\\
& \leq-\phi {} U_{t t}^r-\phi {}\left(U_{t t}-\bar{\lambda}^2 U_{x x}\right)-2 \phi_t {} U_{t t}^r-2 \phi_t {} U_{t t}+2 \bar{\lambda}^2 \phi_t {} U_{x x}+\left(\phi+2 \phi_t\right) {} F.
\end{align}
\noindent Multiplying (\ref{div form: degenerate})   by $(1+t)^\gamma$ ($\gamma>0$ large enough is to be chosen), and then integrating it over $[0, t]\times \mathbb{R}$, one could get
\begin{align}
\label{initial weighted esti: degenertae case}
& (1+t)^\gamma\left\|\left(\phi, \phi_x, \phi_t\right)(t)\right\|_{L^2}^2+\int_0^t(1+\tau)^\gamma\left(\left\|\phi_x(\tau)\right\|_{L^2}^2+\left\|\phi_t(\tau)\right\|_{L^2}^2+G(\tau)\right) d \tau \notag\\
\leq & C\left(\left\|\phi_0\right\|_{H^1}^2+\left\|\phi_1\right\|_{L^2}^2\right)+C \int_0^t(1+\tau)^{\gamma-1}\|\phi(\tau)\|_{L^2}^2 d \tau+C \int_0^t(1+\tau)^{\gamma-1}\left(\left\|\phi_t(\tau)\right\|_{L^2}^2+\left\|\phi_x(\tau)\right\|_{L^2}^2\right) d \tau \notag\\
+ & \int_0^t(1+\tau)^\gamma\left|\int_{\mathbb{R}} \phi_t\left(g^{\prime}(\widetilde{U}+\phi)-g^{\prime}(\widetilde{U})\right) \widetilde{U}_x d x\right| d \tau+\int_0^t(1+\tau)^\gamma\left|\int_{\mathbb{R}} \phi\left\{-\left(U_{t t}-\bar{\lambda}^2 U_{x x}\right)+F-U_{t t}^r\right\} d x\right| d \tau \notag\\
+ & \int_0^t(1+\tau)^\gamma \int_{\mathbb{R}}\left|\phi_t\right|\left|-2 U_{t t}^r-2 U_{t t}+2 \bar{\lambda}^2 U_{x x}+2 F\right| d x d \tau=: C\left(\left\|\phi_0\right\|_{L^2}^2+\left\|\phi_1\right\|_{L^2}^2+\left\|\phi_{0 x}\right\|_{L^2}^2\right)+\sum_{i=1}^5 J_i.
\end{align}
\noindent where we have used the fact that 
\begin{equation}
\int_{\mathbb{R}}\left(g(\widetilde{U}+\phi)-g(\widetilde{U})-g^{\prime}(\widetilde{U}) \phi\right) \widetilde{U}_x(\tau, x) d x \geq C {} G(\tau),
\end{equation}
\noindent which can be seen in \cite{MR4983541}. Now, we estimate $J_1-J_5$ term by term. At first, using the same argument in (\ref{I1 esti}), (\ref{I2 esti}) and (\ref{I5 esti}),  it follows from Lemma \ref{Properties of the smooth rarefaction wave 2}, Lemma \ref{Properties of the viscous contact wave}, Lemma \ref{uniform esti: degenerate case}, Lemma \ref{L1 esti: degenerate case} and Lemma \ref{F(U,U^r) properties} that for $\gamma>2$,
\begin{equation}
    \label{J1 esti}
J_1 \leq \varepsilon \int_0^t(1+\tau)^\gamma\left\|\phi_x(\tau)\right\|_{L^2}^2 d \tau+C(1+t)^{\gamma+\delta-\frac{1}{2}},
\end{equation}

\begin{equation}
    \label{J2 esti}
J_2 \leq C(1+t)^{\gamma-1} \int_0^t\left(\left\|\phi_t\right\|_{L^2}^2+\left\|\phi_x\right\|_{L^2}^2\right) d \tau \leq C(1+t)^{\gamma-1},
\end{equation}
\begin{equation}
    \label{J5 esti}
\begin{aligned}
J_5 & \leq \varepsilon \int_0^t(1+\tau)^\gamma\left\|\phi_t(\tau)\right\|_{L^2}^2 d \tau+C \int_0^t(1+\tau)^\gamma\left(\|F\|{ }_{L^2}^2+\left\|U_{t t}^r\right\|_{L^2}^2+\left\|U_{t t}\right\|_{L^2}^2+\left\|U_{x x}\right\|_{L^2}^2\right) d \tau \\
& \leq \varepsilon \int_0^t(1+\tau)^\gamma\left\|\phi_t(\tau)\right\|_{L^2}^2 d \tau+C \int_0^t(1+\tau)^{\gamma-\frac{3}{2}+\delta} d \tau \\
& \leq \varepsilon \int_0^t(1+\tau)^\gamma\left\|\phi_t(\tau)\right\|_{L^2}^2 d \tau+C(1+t)^{\gamma-\frac{1}{2}+\delta}.
\end{aligned}
\end{equation}
\noindent Now we turn to estimate $J_4$. 
\begin{equation}
\label{J4 esti: initial}
\begin{aligned}
J_4 \leq & \int_0^t(1+\tau)^\gamma \int_{\mathbb{R}}|\phi| {}\left|U_{t t}-\bar{\lambda}^2 U_{x x}\right|(\tau, x) d x d \tau+\int_0^t(1+\tau)^\gamma \int_{\mathbb{R}}|\phi| {}\left(|F|+\left|U_{t t}^r\right|\right)(\tau, x) d x d \tau .
\end{aligned}
\end{equation}
\noindent Denote that 
$W^C(\tau,y)=U(1+\tau, y+\bar{\lambda}(1+\tau))$, $\Phi(\tau,y)=\phi(\tau, y+\bar{\lambda}(1+\tau))$, then under the changing variables
 that $y=x-\bar{\lambda}(1+ t)$, $\tau=t$, we have
\begin{equation*}
\begin{aligned}
    U_{t t}(1+t, x) & =\partial_t\left(\partial_\tau U(1+\tau, y+\bar{\lambda}(1+\tau))-\bar{\lambda} \partial_y U(1+\tau, y+\bar{\lambda}(1+\tau))\right) \\
& =\bar{\lambda}^2 \partial_{y y} U(1+\tau, y+\bar{\lambda}(1+\tau))-2 \bar{\lambda} \partial_{\tau y} U(1+\tau, y+\bar{\lambda}(1+\tau))+\partial_{\tau \tau} U(1+\tau, y+\bar{\lambda}(1+\tau)) \\
& =\bar{\lambda}^2 W^C_{y y}-2 \bar{\lambda} W^C_{\tau y}+W^C_{\tau \tau},
\end{aligned}
\end{equation*}
\noindent and
\begin{equation*}
    U_{x x}(1+t, x)  =\partial_{y y} U(1+\tau, y+\bar{\lambda}(1+\tau))=W^C_{y y}(\tau, y).
\end{equation*}
\noindent From Lemma \ref{Properties of the viscous contact wave} we know that 
$$
\begin{gathered}\left\|W^C_{\tau y}(\tau)\right\|_{L^1}=\left\|\partial_{\tau y} U(1+\tau, y+\bar{\lambda}(1+\tau))(\tau)\right\|_{L^1} \leq C(1+\tau)^{-1} ,\\ \left\|W^C_{\tau \tau}(\tau)\right\|_{L^1}=\left\|\partial_{\tau \tau} U(1+\tau, y+\bar{\lambda}(1+\tau))(\tau)\right\|_{L^1} \leq C(1+\tau)^{-\frac{3}{2}},\end{gathered}
$$
\noindent thus one could obtain that 
\begin{equation}
\label{changing variables esti}
    \begin{aligned}
\int_{\mathbb{R}}\left|U_{t t}-\bar{\lambda}^2 U_{x x}\right|(t, x) d x & =\int_{\mathbb{R}}\left|\bar{\lambda}^2 W_{y y}^C-2 \bar{\lambda} W_{\tau y}^C+W_{\tau \tau}^C-\bar{\lambda}^2 W_{y y}^C\right|(\tau, y) d y \\
& =\int_{\mathbb{R}}\left|-2 \bar{\lambda} W_{\tau y}^C+W_{\tau \tau}^C\right|(\tau, y) d y \leq C\left(\left\|W_{\tau y}^C(\tau)\right\|_{L^1}+\left\|W_{\tau \tau}^C(\tau)\right\|_{L^1}\right) \\
& \leq C(1+\tau)^{-1}=C(1+t)^{-1}.
\end{aligned}
\end{equation}

\noindent Then by employing the same argument as in (\ref{I4 esti}), we can have
\begin{equation}
\label{J4 esti}
\begin{aligned}
J_4 & \leq C \int_0^t(1+\tau)^\gamma\|\phi(\tau)\|_{L^{\infty}} {}\left(\left\|U_{t t}-\bar{\lambda}^2 U_{x x}(\tau)\right\|_{L^1}+\|F(\tau)\|_{L^1}+\left\|U_{t t}^r(\tau)\right\|_{L^1}\right) d \tau \\
& \leq C \int_0^t(1+\tau)^\gamma\|\phi(\tau)\|_{L^{\infty}} {}\left((1+\tau)^{-1}+(1+\tau)^{-1+\frac{\delta}{2}}\right) d \tau \\
& \leq C \int_0^t(1+\tau)^{\gamma-1+\frac{\delta}{2}}\|\phi(\tau)\|_{L^{\infty}} d \tau \\
& \leq C \int_0^t(1+\tau)^{\gamma-1+\frac{\delta}{2}}\|\phi(\tau)\|_{L^1}^{1 / 3} {}\left\|\phi_x(\tau)\right\|_{L^2}^{2 / 3} d \tau \\
& \leq C \int_0^t(1+\tau)^{\gamma-1+\frac{5}{6} \delta}\left\|\phi_x(\tau)\right\|_{L^2}^{2 / 3} d \tau \\
& \leq \varepsilon \int_0^t(1+\tau)^\gamma\left\|\phi_x(\tau)\right\|_{L^2}^2 d \tau+C(1+t)^{\gamma-\frac{1}{2}+\frac{5}{4} \delta}.
\end{aligned}
\end{equation}
\noindent where we have also used Lemma \ref{F(U,U^r) properties}, Lemma \ref{Properties of the smooth rarefaction wave 2}, and the Young inequality. Subsequently, owing to (6) of Lemma \ref{properties of $g$}, we can deduce that
    \begin{align}
        \label{J3 esti}
J_3 & \leq \varepsilon \int_0^t(1+\tau)^\gamma\left\|\phi_t\right\|_{L^2}^2 d \tau+C \int_0^t(1+\tau)^\gamma \int_{\mathbb{R}}\left(g^{\prime}(\widetilde{U}+\phi)-g^{\prime}(\widetilde{U})\right)^2 {} \widetilde{U}_x^2(\tau, x) d x d \tau \notag\\
& \leq \varepsilon \int_0^t(1+\tau)^\gamma\left\|\phi_t\right\|_{L^2}^2 d \tau+C \int_0^t(1+\tau)^\gamma\left\|\widetilde{U}_x\right\|_{L^{\infty}}^2 \int_{\mathbb{R}}\left(g^{\prime}(\widetilde{U}+\phi)-g^{\prime}(\widetilde{U})\right)^2(\tau, x) d x d \tau \notag\\
& \leq \varepsilon \int_0^t(1+\tau)^\gamma\left\|\phi_t\right\|_{L^2}^2 d \tau+C \int_0^t(1+\tau)^{\gamma-1} \int_{\mathbb{R}}\left(g^{\prime}(\widetilde{U}+\phi)-g^{\prime}(\widetilde{U})\right)^2(\tau, x) d x d \tau \notag\\
&\leq \varepsilon \int_0^t(1+\tau)^\gamma\left\|\phi_t\right\|_{L^2}^2 d \tau+C \int_0^t(1+\tau)^{\gamma-1}\|\phi(\tau)\|_{L^2}^2 d \tau \notag\\
&\leq  \varepsilon \int_0^t(1+\tau)^\gamma\left\|\phi_t\right\|_{L^2}^2 d \tau+J_1\notag\\
&\leq  \varepsilon \int_0^t(1+\tau)^\gamma\left(\left\|\phi_t\right\|_{L^2}^2+\left\|\phi_x\right\|_{L^2}^2\right)d \tau+C(1+t)^{\gamma+\delta-1 / 2},
\end{align}
\noindent Inserting (\ref{J1 esti}), (\ref{J2 esti}), (\ref{J3 esti}), (\ref{J4 esti}), (\ref{J5 esti}) into (\ref{initial weighted esti: degenertae case}) and take $\varepsilon>0$ sufficiently small, it holds that
\begin{equation}
    \begin{aligned}
& (1+t)^\gamma\left(\|\phi(t)\|_{L^2}^2+\left\|\phi_t(t)\right\|_{L^2}^2+\left\|\phi_x(t)\right\|_{L^2}^2\right)+\int_0^t(1+\tau)^\gamma\left(\left\|\phi_t(\tau)\right\|_{L^2}^2+\left\|\phi_x(\tau)\right\|_{L^2}^2+G(\tau)\right) d \tau \\
& \leq C\left(\left\|\phi_0\right\|_{L^2}^2+\left\|\phi_{0 x}\right\|_{L^2}^2+\left\|\phi_1\right\|_{L^2}^2\right)+C(1+t)^{\gamma-\frac{1}{2}+\frac{5}{4} \delta}.
\end{aligned}
\end{equation}

\end{proof}

\subsection{$L^2$ decay for $(\phi_x, \psi_x)$}
\begin{lemma}
    \label{derivative L2 decay: degenerate case}
Under the same assumptions as in Theorem \ref{Main Th 2.1}, for $\delta>0$  as in Lemma \ref{F(U,U^r) properties} and $\gamma>4$,
there exists a constant $C>0$ 
such that for $t\geq 0$,
\begin{equation*}
    \begin{aligned}
&(1+t)^\gamma\left(\left\|\phi_x(t)\right\|_{L^2}^2+\left\|\phi_{t x}(t)\right\|_{L^2}^2+\left\|\phi_{x x}(t)\right\|_{L^2}^2\right)+\int_0^t(1+\tau)^\gamma\left(\left\|\phi_{t x}(\tau)\right\|_{L^2}^2+\left\|\phi_{x x}(\tau)\right\|_{L^2}^2\right) d \tau \\
& \leq C(\phi_0, \phi_1)+C(1+t)^{\gamma-\frac{5}{6}+2\delta}.
\end{aligned}
\end{equation*}
\end{lemma}
\noindent The proof of Lemma \ref{derivative L2 decay: degenerate case} is technically involved and constitutes a main novelty of this paper.  Before the proof of Lemma \ref{derivative L2 decay: degenerate case}, we introduce the Weak derivative formula for composite functions. It plays an essential role in dealing with the partially linearly degenerate flux.
\begin{lemma}
\label{composite derivative}
 Given $f\in C^0(\mathbb{R})$, and $f$ is continuously differentiable in $\mathbb{R}\backslash L$, where $L$ is a discrete subset of $\mathbb{R}$. Assume also that $f'$ is bounded in $\mathbb{R}\backslash L$. Let $\Omega\subset\mathbb R^n$, $n\geq1$, be open and $u\in W^{1,p}_{\mathrm{loc}}(\Omega)$ for $1\leq p<\infty$. Then $f\circ u\in W^{1,p}_{\mathrm{loc}}(\Omega)$ and, for $1\leq i\leq n$, the following formula holds almost everywhere:
 \begin{equation*}
     \partial_i(f \circ u)(x)= \begin{cases}f^{\prime}(u(x)) \partial_i u(x), & u(x) \notin L, \\ 0, & u(x) \in L .\end{cases}
 \end{equation*}
\end{lemma}
\begin{proof}
The assumptions imply that $f$ is globally Lipschitz. Apply the Sobolev chain rule for Lipschitz functions. For each $c\in L$, $\partial_i u=0$ almost everywhere on $\{u=c\}$. Since $L$ is countable, this proves the stated formula also on $\{u\in L\}$.
\end{proof}
\begin{remark}
    Considering a special case, if $\|u\|_{L^{\infty}(\Omega)} \leq C<+\infty$ and $L=\{0\}$, the conclusion still holds if $f'$ is bounded on $[-C,C]\backslash L$.
\end{remark}
\begin{proof}
Define
\begin{equation*}
\tilde f(v)=
\begin{cases}
f(-C), & v<-C,\\
f(v), & |v|\leq C,\\
f(C), & v>C.
\end{cases}
\end{equation*}
Then $\tilde f$ is globally Lipschitz, its derivative is bounded away from $\{-C,0,C\}$, and $\tilde f(u)=f(u)$. Applying Lemma \ref{composite derivative} to $\tilde f$ with $L$ replaced by $\{-C,0,C\}$ gives, for a.e. $x\in\Omega$,
\begin{equation}
    \begin{aligned}
\partial_{x_i} f(u(x))=\partial_{x_i} \tilde{f}(u(x)) & = \begin{cases}\tilde{f}^{\prime}(u(x)) \partial_i u(x), & u \neq 0, \\
0, & u=0,\end{cases} \\
& = \begin{cases}f^{\prime}(u(x)) \partial_i u(x), & u \neq 0, \\
0, & u=0.\end{cases}
\end{aligned}
\end{equation}
\end{proof}

\noindent \textbf{Proof of Lemma \ref{derivative L2 decay: degenerate case}.}

In the end, we will derive the refined decay rate of $\left\|\phi_x(t)\right\|_{L^2}$ by using the time-weighted $L^2$ method. Multiplying (\ref{perturbation equ: degenerate case}) by $-\phi_{xx}$,
\begin{equation}
    \label{mutiply -phi_xx: degenerate case}
    \begin{aligned}
& \left\{\phi_x \phi_{t x}+\frac{1}{2} \phi_x^2\right\}_t-\phi_{x x} {}\left(g^{\prime}(\widetilde{U}+\phi)-g^{\prime}(\widetilde{U})\right) \widetilde{U}_x-\phi_{x x} {}\left(g^{\prime}(\widetilde{U}+\phi)-g^{\prime}(\widetilde{U})\right) \phi_x \\
& \quad+\left\{-\phi_x \phi_{t t}-\phi_x \phi_t-\frac{1}{2} \phi_x^2 {} g^{\prime}(\widetilde{U})-\frac{\bar{\lambda}}{2} \phi_x^2\right\}_x+\phi_{x x}^2-\phi_{t x}^2 \leq \phi_{x x} U_{t t}^r+\phi_{x x} {}\left(U_{t t}-\bar{\lambda}^2 U_{x x}\right)-\phi_{x x} F.
\end{aligned}
\end{equation}
\noindent where we have used the fact that
    \begin{align}
-\phi_{x x} g^{\prime}(\widetilde{U}+\phi) \phi_x & =-\phi_{x x} {}\left(g^{\prime}(\widetilde{U}+\phi)-g^{\prime}(\widetilde{U})\right) \phi_x-\phi_{x x} g^{\prime}(\widetilde{U}) \phi_x \notag\\
& =-\phi_{x x} {}\left(g^{\prime}(\widetilde{U}+\phi)-g^{\prime}(\widetilde{U})\right) \phi_x+\left\{-\frac{1}{2} \phi_x^2 g^{\prime}(\widetilde{U})\right\}_x+\frac{1}{2} \phi_x^2 g^{\prime \prime}(\widetilde{U}) \widetilde{U}_x \notag\\
& \geq-\phi_{x x} {}\left(g^{\prime}(\widetilde{U}+\phi)-g^{\prime}(\widetilde{U})\right) \phi_x+\left\{-\frac{1}{2} \phi_x^2 g^{\prime}(\widetilde{U})\right\}_x.
\end{align}
\noindent Then, for fixed $t>0$,
denote $\Omega_0:=\{x \in \mathbb{R} \mid(\widetilde{U}+\phi)(t, x)=0\}, \Omega_1:=\{x \in \mathbb{R} \mid(\widetilde{U}+\phi)(t, x) \neq 0\}$, it is obvious that $\Omega_0 \cup \Omega_1=\mathbb{R}$. Noticing that $g^\prime(u)\in C^0(\mathbb{R})$ and $g^{\prime\prime}(u)$ is bounded on $[-(B(N_0)+\max( u_+, -u_-)), B(N_0)+\max( u_+, -u_-)]\backslash \{0\}$, owing to Lemma \ref{composite derivative}, a direct calculation gives
    \begin{align}
& (g(\widetilde{U}+\phi)-g(\widetilde{U}))_{x x} \notag\\
= & \left(g^{\prime \prime}(\widetilde{U}+\phi)-g^{\prime \prime}(\widetilde{U})\right) \widetilde{U}_x^2+2 g^{\prime \prime}(\widetilde{U}+\phi) \phi_x \widetilde{U}_x+\left(g^{\prime}(\widetilde{U}+\phi)-g^{\prime}(\widetilde{U})\right) \widetilde{U}_{x x} \notag\\
& +g^{\prime}(\widetilde{U}+\phi) \phi_{x x}+g^{\prime \prime}(\widetilde{U}+\phi) \phi_x^2, \quad \mathrm{a.e.} x \in \Omega_1, 
\end{align}
\noindent and 
    \begin{align}
& (g(\widetilde{U}+\phi)-g(\widetilde{U}))_{x x} \notag\\
= & -g^{\prime \prime}(\widetilde{U}) \widetilde{U}_x^2-g^{\prime}(\widetilde{U}) \widetilde{U}_{x x}+g^{\prime}(\widetilde{U}+\phi)\left(\widetilde{U}_{x x}+\phi_{x x}\right)+\left(g^{\prime}(\widetilde{U}+\phi)\right)_x\left(\widetilde{U}_x+\phi_x\right) \notag\\
= & -g^{\prime \prime}(\widetilde{U}) \widetilde{U}_x^2-g^{\prime}(\widetilde{U}) \widetilde{U}_{x x}+0 \cdot{}\left(\widetilde{U}_{x x}+\phi_{x x}\right)+0\cdot {}\left(\widetilde{U}_x+\phi_x\right) \notag\\
= & -g^{\prime \prime}(\widetilde{U}) \widetilde{U}_x^2-g^{\prime}(\widetilde{U}) \widetilde{U}_{x x}, \quad \mathrm{a.e.} x \in \Omega_0.
\end{align}
\noindent Thus, taking weak derivative on (\ref{perturbation equ: degenerate case}) with respect to $x$, we have 
\begin{equation}
\label{derivative1: degenerate case}
    \begin{aligned}
& \phi_{t t x}+\phi_{t x}+\left(g^{\prime \prime}(\widetilde{U}+\phi)-g^{\prime \prime}(\widetilde{U})\right) \widetilde{U}_x^2+2 g^{\prime \prime}(\widetilde{U}+\phi) \phi_x \widetilde{U}_x+\left(g^{\prime}(\widetilde{U}+\phi)-g^{\prime}(\widetilde{U})\right) \widetilde{U}_{x x} \\
& +g^{\prime}(\widetilde{U}+\phi) \phi_{x x}+g^{\prime \prime}(\widetilde{U}+\phi) \phi_x^2+\bar{\lambda} \phi_{x x}-\phi_{x x x}=\partial_x F\left(U, U^r\right)-U_{t t x}^r-U_{t t x}+\bar{\lambda}^2 U_{x x x}, \quad \mathrm{a.e.} x \in \Omega_1
\end{aligned}
\end{equation}
\noindent and
\begin{equation}
\label{derivative 2: degenerate case}
    \phi_{t t x}+\phi_{t x}-g^{\prime \prime}(\widetilde{U}) \widetilde{U}_x^2-g^{\prime}(\widetilde{U}) \widetilde{U}_{x x}+\bar{\lambda} \phi_{x x}-\phi_{x x x}=\partial_x F\left(U, U^r\right)-U_{t t x}^r-U_{t t x}+\bar{\lambda}^2 U_{x x x},  \quad \text { a.e.} x \in \Omega_0.
\end{equation}
\noindent Multiplying (\ref{derivative1: degenerate case}) and (\ref{derivative 2: degenerate case}) by $2\phi_{tx}$ respectively, we deduce that
\begin{equation}
    \label{derivative 11: degenertae case}
\begin{aligned}
& \left(\phi_{t x}^2+\phi_{x x}^2\right)_t+\left\{-2 \phi_{t x} \phi_{x x}\right\}_x+2 \phi_{t x}^2+2\left(g^{\prime}(u)+\bar{\lambda}\right) \phi_{t x} {} \phi_{x x} \\
+ & 2 \phi_{t x} {}\left(\left(g^{\prime \prime}(\widetilde{U}+\phi)-g^{\prime \prime}(\widetilde{U})\right) \widetilde{U}_x^2+2 g^{\prime \prime}(\widetilde{U}+\phi) \phi_x \widetilde{U}_x+\left(g^{\prime}(\widetilde{U}+\phi)-g^{\prime}(\widetilde{U})\right) \widetilde{U}_{x x}+g^{\prime \prime}(\widetilde{U}+\phi) \phi_x^2\right) \\
= & 2 \phi_{t x} {} \partial_x F\left(U, U^r\right)-2 \phi_{t x} {} U_{t t x}^r-2 \phi_{t x} {} U_{t t x}+2 \bar{\lambda}^2 \phi_{t x} {} U_{x x x}, \quad \text { a.e. } x \in \Omega_1
\end{aligned}
\end{equation}
\noindent and
\begin{equation}
    \label{derivative 21: degenerate case}
\begin{aligned}
& \left(\phi_{t x}^2+\phi_{x x}^2\right)_t+\left\{-2 \phi_{t x} \phi_{x x}\right\}_x+2 \phi_{t x}^2+2\left(g^{\prime}(0)+\bar{\lambda}\right) \phi_{t x} {} \phi_{x x}+2 \phi_{t x} {}\left(-g^{\prime \prime}(\widetilde{U}) \widetilde{U}_x^2-g^{\prime}(\widetilde{U}) \widetilde{U}_{x x}\right) \\
& =2 \phi_{t x} {} \partial_x F\left(U, U^r\right)-2 \phi_{t x} {} U_{t t x}^r-2 \phi_{t x} {} U_{t t x}+2 \bar{\lambda}^2 \phi_{t x} {} U_{x x x}, \quad \text { a.e. } x \in \Omega_0.
\end{aligned}
\end{equation}
\noindent First combine (\ref{derivative 11: degenertae case}) and (\ref{derivative 21: degenerate case}) into an almost-everywhere identity on $\mathbb R$, and add (\ref{mutiply -phi_xx: degenerate case}). Then multiply by $(1+t)^\gamma$ and integrate over $[0,t]\times\mathbb R$. The resulting estimate is
\begin{align}
\label{derivative weighted esti initial: degenerate case}
& (1+t)^\gamma\left(\left\|\phi_x(t)\right\|_{L^2}^2+\left\|\phi_{t x}(t)\right\|_{L^2}^2+\left\|\phi_{x x}(t)\right\|_{L^2}^2\right)+\int_{0}^t(1+\tau)^\gamma\left(\left\|\phi_{x x}(\tau)\right\|_{L^2}^2+\left\|\phi_{t x}(\tau)\right\|_{L^2}^2\right) d \tau \notag\\
\leq & C(\phi_0,\phi_1)+C \int_{0}^t(1+\tau)^{\gamma-1}\left(\left\|\phi_x(\tau)\right\|_{L^2}^2+\left\|\phi_{t x}(\tau)\right\|_{L^2}^2+\left\|\phi_{x x}(\tau)\right\|_{L^2}^2\right) d \tau \notag\\
 +& C \int_{0}^t(1+\tau)^\gamma \int_{\mathbb{R}}\left|\phi_{x x}\left(g^{\prime}(\widetilde{U}+\phi)-g^{\prime}(\widetilde{U})\right)\right| \widetilde{U}_x d x d \tau+C \int_{0}^t(1+\tau)^\gamma \int_{\mathbb{R}}\left|\phi_{x x}\left(g^{\prime}(\widetilde{U}+\phi)-g^{\prime}(\widetilde{U})\right) \phi_x\right| d x d \tau \notag\\
+& C \int_{0}^t(1+\tau)^\gamma \int_{\Omega_1}\left|\phi_{t x}\left(g^{\prime \prime}(\widetilde{U}+\phi)-g^{\prime \prime}(\widetilde{U})\right)\right| \widetilde{U}_x^2 d x d \tau+C \int_{0}^t(1+\tau)^\gamma \int_{\Omega_1} \mid \phi_{t x} g^{\prime \prime}(\widetilde{U}+\phi) \phi_x\widetilde{U}_x\mid  d x d \tau \notag\\
+& C \int_{0}^t(1+\tau)^\gamma \int_{\Omega_1}\left|\phi_{t x}\left(g^{\prime}(\widetilde{U}+\phi)-g^{\prime}(\widetilde{U})\right) \widetilde{U}_{x x}\right| d x d \tau+C \int_{0}^t(1+\tau)^\gamma \int_{\Omega_1}\left|\phi_{t x} g^{\prime \prime}(\widetilde{U}+\phi) \phi_x^2\right| d x d \tau \notag\\
+& C \int_{0}^t(1+\tau)^\gamma \int_{\Omega_0}\left|\phi_{t x}\left(-g^{\prime \prime}(\widetilde{U}) \widetilde{U}_x^2-g^{\prime}(\widetilde{U}) \widetilde{U}_{x x}\right)\right| d x d \tau \notag\\
+& C \int_{0}^t(1+\tau)^\gamma \int_{\mathbb{R}}\left|\phi_x F_x\right| d x d \tau+C \int_{0}^t(1+\tau)^\gamma \int_{\mathbb{R}}\left|\phi_{x x}\left(U_{t t}-\bar{\lambda}^2 U_{x x}+U_{t t}^r\right)\right| d x d \tau \notag\\
+& C \int_{0}^t(1+\tau)^\gamma \int_{\mathbb{R}}\left|\phi_{t x}\right| {}\left(\left|\partial_x F\left(U, U^r\right)\right|+\left|U_{t t x}^r\right|+\left|U_{t t x}\right|+\left|U_{x x x}\right|\right) d x d \tau \notag\\
& =: C(\phi_0,\phi_1)+\sum_{i=6}^{16} J_i.
\end{align}
\noindent  The uniform energy estimate in Lemma \ref{uniform esti: degenerate case} gives
\begin{equation}
    \phi_{t x}^2+\phi_{x x}^2+2\left(\bar{\lambda}+g^{\prime}(u)\right) \phi_{t x} \phi_{x x} \geq \alpha\left(\phi_{t x}^2+\phi_{x x}^2\right)(\alpha>0),
\end{equation}
\noindent which can be seen in \cite{MR4983541}.  Next  we will estimate $J_6-J_{16}$ term by term. By applying the same argument as in (\ref{I6 esti}), we have
\begin{equation}
\label{J6 esti}
    J_6 \leq \varepsilon \int_{0}^t(1+\tau)^\gamma\left(\left\|\phi_{t x}\right\|_{L^2}^2+\left\|\phi_{x x}\right\|_{L^2}^2\right) d \tau+C(1+t)^{\gamma-\frac{3}{2}+\frac{5}{4}\delta}.
\end{equation}
\noindent Subsequently, applying the same estimate as in (\ref{J3 esti}), one could have  
    \begin{align}
    \label{J8,J11,J12 esti}
& J_8+J_{11}+J_{12} \notag\\
\leq & \varepsilon \int_{0}^t(1+\tau)^\gamma\left(\left\|\phi_{x x}(\tau)\right\|_{L^2}^2+\left\|\phi_{t x}(\tau)\right\|_{L^2}^2\right) d \tau+C \int_{0}^t(1+\tau)^\gamma \int_{\mathbb{R}}\left(g^{\prime}(\widetilde{U}+\phi)-g^{\prime}(\widetilde{U})\right)^2 {} \phi_x^2 d x d \tau \notag\\
& +C \int_{0}^t(1+\tau)^\gamma \int_{\mathbb{R}}\left(g^{\prime}(\widetilde{U}+\phi)-g^{\prime}(\widetilde{U})\right)^2 {} \widetilde{U}_{x x}^2 d x d \tau+\int_{0}^t(1+\tau)^\gamma \int_{\mathbb{R}} \phi_x^4 d x d \tau \notag\\
\leq & \varepsilon \int_{0}^t(1+\tau)^\gamma\left(\left\|\phi_{x x}(\tau)\right\|_{L^2}^2+\left\|\phi_{t x}(\tau)\right\|_{L^2}^2\right) d \tau+C \int_{0}^t(1+\tau)^\gamma \int_{\mathbb{R}} \phi^2 {} \phi_x^2 d x d \tau \notag\\
& +C \int_{0}^t(1+\tau)^\gamma\left\|\widetilde{U}_{x x}(\tau)\right\|_{L^{\infty}}^2\|\phi(\tau)\|_{L^2}^2 d \tau+\int_{0}^t(1+\tau)^\gamma \int_{\mathbb{R}} \phi_x^4 d x d \tau \notag\\
\leq & \varepsilon \int_{0}^t(1+\tau)^\gamma\left(\left\|\phi_{x x}(\tau)\right\|_{L^2}^2+\left\|\phi_{t x}(\tau)\right\|_{L^2}^2\right) d \tau+C(1+t)^{\gamma-\frac{3}{2}+7 \delta}+C(1+t)^{\gamma-\frac{3}{2}+4 \delta}+C(1+t)^{\gamma-\frac{3}{2}+\frac{15}{4} \delta} \notag\\
\leq & \varepsilon \int_{0}^t(1+\tau)^\gamma\left(\left\|\phi_{x x}(\tau)\right\|_{L^2}^2+\left\|\phi_{t x}(\tau)\right\|_{L^2}^2\right) d \tau+C(1+t)^{\gamma-\frac{3}{2}+7 \delta}.
\end{align}
\noindent Similarly, it holds from Lemma \ref{L2 decay: degenerate case} that for $\gamma>4$, 
\begin{align}
    \label{J7 esti}
    J_7 & \leq \varepsilon \int_{0}^t(1+\tau)^\gamma\left\|\phi_{x x}(\tau)\right\|_{L^2}^2 d \tau+C \int_{0}^t(1+\tau)^\gamma \int_{\mathbb{R}}\left(g^{\prime}(\widetilde{U}+\phi)-g^{\prime}(\widetilde{U})\right)^2 {} \widetilde{U}_x^2 d x d \tau \notag\\
& \leq \varepsilon \int_{0}^t(1+\tau)^\gamma\left\|\phi_{x x}(\tau)\right\|_{L^2}^2 d \tau+C \int_{0}^t(1+\tau)^{\gamma-\frac{1}{2}} \int_{\mathbb{R}}\left(g^{\prime}(\widetilde{U}+\phi)-g^{\prime}(\widetilde{U})\right)^2 {} \widetilde{U}_x d x d \tau \notag\\
& \leq \varepsilon \int_{0}^t(1+\tau)^\gamma\left\|\phi_{x x}(\tau)\right\|_{L^2}^2 d \tau+C \int_{0}^t(1+\tau)^{\gamma-\frac{1}{2}} G(\tau) d \tau \notag\\
& =\varepsilon \int_{0}^t(1+\tau)^\gamma\left\|\phi_{x x}(\tau)\right\|_{L^2}^2 d \tau+C \int_{0}^t(1+\tau)^{\gamma-3-\frac{1}{2}} {}(1+\tau)^3 G(\tau) d \tau \notag\\
& \leq \varepsilon \int_{0}^t(1+\tau)^\gamma\left\|\phi_{x x}(\tau)\right\|_{L^2}^2 d \tau+C(1+t)^{\gamma-3-\frac{1}{2}} \int_{0}^t(1+\tau)^3 G(\tau) d \tau \notag\\
& \leq \varepsilon \int_{0}^t(1+\tau)^\gamma\left\|\phi_{x x}(\tau)\right\|_{L^2}^2 d \tau+C(1+t)^{\gamma-3-\frac{1}{2}}(1+t)^{3-\frac{1}{2}+\frac{5}{4}\delta} \notag\\
& \leq \varepsilon \int_{0}^t(1+\tau)^\gamma\left\|\phi_{x x}(\tau)\right\|_{L^2}^2 d \tau+C(1+t)^{\gamma-1+\frac{5}{4}\delta}.
\end{align}

\noindent Moreover, the Young inequality, Lemma \ref{Properties of the viscous contact wave}, Lemma \ref{F(U,U^r) properties} together imply that
\begin{align}
   \label{J10,J15,J16 esti}
    & J_{10}+J_{15}+J_{16} \notag\\
\leq & \varepsilon \int_{0}^t(1+\tau)^\gamma\left(\left\|\phi_{x x}(\tau)\right\|_{L^2}^2+\left\|\phi_{t x}(\tau)\right\|_{L^2}^2\right) d \tau+C \int_{0}^t(1+\tau)^\gamma \int_{\mathbb{R}} \phi_x^2 {} \widetilde{U}_x^2 d x d \tau \notag\\
& +C \int_{0}^t(1+\tau)^\gamma\left(\left\|\partial_x F\right\|_{L^2}^2+\left\|U_{t t x}^r\right\|_{L^2}^2+\left\|U_{t t x}\right\|_{L^2}^2+\left\|U_{x x x}\right\|_{L^2}^2\right) d \tau \notag\\
& +C \int_{0}^t(1+\tau)^\gamma\left(\left\|U_{t t}^r\right\|_{L^2}^2+\left\|U_{t t}-\bar{\lambda}^2 U_{x x}(\tau)\right\|_{L^2}^2\right) d \tau \notag\\
\leq & \varepsilon \int_{0}^t(1+\tau)^\gamma\left(\left\|\phi_{x x}(\tau)\right\|_{L^2}^2+\left\|\phi_{t x}(\tau)\right\|_{L^2}^2\right) d \tau+\int_{0}^t(1+\tau)^{\gamma-1}\left\|\phi_x(\tau)\right\|_{L^2}^2 d \tau+C(1+t)^{\gamma-1+\delta} \notag\\
\leq & \varepsilon \int_{0}^t(1+\tau)^\gamma\left(\left\|\phi_{x x}(\tau)\right\|_{L^2}^2+\left\|\phi_{t x}(\tau)\right\|_{L^2}^2\right) d \tau+J_6+C(1+t)^{\gamma-1+\delta} \notag\\
\leq & \varepsilon \int_{0}^t(1+\tau)^\gamma\left(\left\|\phi_{x x}(\tau)\right\|_{L^2}^2+\left\|\phi_{t x}(\tau)\right\|_{L^2}^2\right) d \tau+C(1+t)^{\gamma-1+\delta},
\end{align}

\noindent and 
        \begin{align}
\label{J13 esti}
J_{13}\leq{}&\varepsilon\int_{0}^t(1+\tau)^\gamma\|\phi_{tx}(\tau)\|_{L^2}^2\,d\tau\notag\\
&+C\int_{0}^t(1+\tau)^\gamma
\left\{\int_{X(\tau)}^\infty\widetilde U_x^4\,dx
+\|g'(\widetilde U)\widetilde U_{xx}\|_{L^2}^2\right\}\,d\tau\notag\\
\leq{}&\varepsilon\int_{0}^t(1+\tau)^\gamma\|\phi_{tx}(\tau)\|_{L^2}^2\,d\tau
+C(1+t)^{\gamma-1+\frac{\delta}{2}}.
\end{align}
\noindent To justify the last step, Lemma \ref{Properties of the smooth rarefaction wave 2}, Lemma \ref{X(t) properties}(5), and (\ref{one residual esti}) give
\begin{align}
\int_{X(t)}^\infty\widetilde U_x^4\,dx
&\leq C\|U_x^r\|_{L^\infty}^3\|U_x^r\|_{L^1}
+C\|U_x\|_{L^\infty}^3[-U(1+t,X(t))]\notag\\
&\leq C(1+t)^{-3}+C(1+t)^{-2+\frac{\delta}{2}},\notag\\
\|g'(\widetilde U)\widetilde U_{xx}\|_{L^2}^2
&\leq\|g'(\widetilde U)\widetilde U_{xx}\|_{L^\infty}
       \|g'(\widetilde U)\widetilde U_{xx}\|_{L^1}\notag\\
&\leq C(1+t)^{-1}
\bigl(\|U_{xx}^r\|_{L^1}+\|g'(\widetilde U)U_{xx}\|_{L^1}\bigr)\notag\\
&\leq C(1+t)^{-2+\frac{\delta}{2}}.
\end{align}
\noindent For the estimate on $J_9$,
\begin{align}
        \label{J9 esti}
J_9 \leq & \varepsilon \int_0^t(1+\tau)^\gamma\left\|\phi_{t x}(\tau)\right\|_{L^2}^2 d \tau+C \int_0^t(1+\tau)^\gamma \int_{\tilde{U}+\phi>0, \widetilde{U}>0}\left(g^{\prime \prime}(\widetilde{U}+\phi)-g^{\prime \prime}(\widetilde{U})\right)^2 {} \widetilde{U}_x^4 d x d \tau \notag\\
& +C \int_0^t(1+\tau)^\gamma \int_{\tilde{U}+\phi>0, \widetilde{U}<0}\left(g^{\prime \prime}(\widetilde{U}+\phi)\right)^2 \widetilde{U}_x^4 d x d \tau+C \int_0^t(1+\tau)^\gamma \int_{\tilde{U}+\phi<0, \widetilde{U}>0}\left(g^{\prime \prime}(\widetilde{U})\right)^2 \widetilde{U}_x^4 d x d \tau \notag\\
\leq & \varepsilon \int_0^t(1+\tau)^\gamma\left\|\phi_{t x}(\tau)\right\|_{L^2}^2 d \tau+C \int_0^t(1+\tau)^\gamma\left\|\widetilde{U}_x(\tau)\right\|_{L^\infty}^4\|\phi(\tau)\|_{L^2}^2 d \tau \notag\\
& +C \int_0^t(1+\tau)^{\gamma-\frac{3}{2}}\|\phi(\tau)\|_{L^\infty} d \tau+C \int_0^t(1+\tau)^\gamma \int_{\tilde{U}>0} \widetilde{U}_x^4 d x d \tau \notag\\
\leq & \varepsilon \int_0^t(1+\tau)^\gamma\left\|\phi_{t x}(\tau)\right\|_{L^{2}}^2 d \tau+C \int_0^t(1+\tau)^{\gamma-2}\|\phi(\tau)\|_{L^2}^2 d \tau+C(1+t)^{\gamma-\frac{13}{14}+\delta}+C(1+t)^{\gamma-1+\frac{\delta}{2}} \notag\\
\leq & \varepsilon \int_0^t(1+\tau)^\gamma\left\|\phi_{t x}(\tau)\right\|_{L^2}^2 d \tau+C(1+t)^{\gamma-\frac{3}{2}+\frac{5}{4} \delta}+C(1+t)^{\gamma-1+\delta},
\end{align}
\noindent where we have used
\begin{align}
\int_{\widetilde{U}+\phi>0, \widetilde{U}<0} \widetilde{U}_x(t, x)^4 d x & \leq\left\|\widetilde{U}_x(t)\right\|_{L^{\infty}}^3 \int_{\widetilde{U}+\phi>0, \widetilde{U}<0} \widetilde{U}_x(t, x) d x \notag\\
& =\left\|\widetilde{U}_x(t)\right\|_{L^{\infty}}^3 \int_{\{\widetilde{U}+\phi>0, \widetilde{U}<0\} \cap\left\{-\phi(t, x) \geq-\|\phi\|_{L^{\infty}}\right\}} \widetilde{U}_x d x \notag\\
& =\left\|\widetilde{U}_x(t)\right\|_{L^{\infty}}^3 \int_{\{-\phi<\widetilde{U}<0\} \cap\left\{-\phi(t, x) \geq-\|\phi\|_{L^{\infty}}\right\}} \widetilde{U}_x d x \notag\\
& \leq\left\|\widetilde{U}_x(t)\right\|_{L^{\infty}}^3 \int_{-\|\phi(t)\|_{L^{\infty}}<\widetilde{U}<0} \widetilde{U}_x d x \notag\\
& \leq C(1+t)^{-\frac{3}{2}}\|\phi(t)\|_{L^{\infty}}.
\end{align}

\noindent and 
    \begin{align}
\int_0^t(1+\tau)^{\gamma-\frac{3}{2}}\|\phi(\tau)\|_{L^{\infty}} d \tau & \leq C \int_0^t(1+\tau)^{\gamma-\frac{3}{2}}\left\|\phi_x(\tau)\right\|{ }_{L^2}^{1 / 2}\|\phi(\tau)\|_{L^2}^{1 / 2} d \tau \notag\\
& \leq C \int_0^t(1+\tau)^{\gamma-\frac{3}{2}}\left\|\phi_{x x}(\tau)\right\|{ }_{L^2}^{1 / 4}\|\phi(\tau)\|_{L^2}^{3 / 4} d \tau \notag\\
& \leq C \int_0^t(1+\tau)^{\gamma-\frac{3}{2}-\frac{3}{16}+\frac{15}{32} \delta}\left\|\phi_{x x}(\tau)\right\|_{L^2}^{1 / 4} d \tau \notag\\
& \leq \varepsilon \int_0^t(1+\tau)^\gamma\left\|\phi_{x x}(\tau)\right\|_{L^2}^2 d \tau+C \int_0^t(1+\tau)^{\gamma-\frac{27}{14}+\delta} d \tau \notag\\
& \leq \varepsilon \int_0^t(1+\tau)^\gamma\left\|\phi_{x x}(\tau)\right\|_{L^2}^2 d \tau+C(1+t)^{\gamma-\frac{13}{14}+\delta}.
\end{align}
\noindent As to $J_{14}$, noting that $\|\partial_xF(t)\|_{L^1}\leq C(1+t)^{-1+\delta}$, Lemma \ref{L2 decay: degenerate case} with $\gamma=3$ yields that 
 \begin{align}
& \int_0^t(1+\tau)^\gamma \int_{\mathbb{R}}\left|\phi_x\right| {}\left|F_x\right| d x d \tau \notag\\
\leq & \varepsilon \int_0^t(1+\tau)^\gamma\left\|\phi_{x x}\right\|_{L^2}^2 d \tau+C \int_0^t(1+\tau)^{\gamma+\frac{3}{2} b} d \tau+C \int_0^t(1+\tau)^{\gamma+3 a}\left\|\phi_x\right\|_{L^2}^2 d \tau \notag\\
\leq & \varepsilon \int_0^t(1+\tau)^\gamma\left\|\phi_{x x}\right\|_{L^2}^2 d \tau+C(1+t)^{\gamma+\frac{3}{2} b+2\delta+1}+C \int_0^t(1+\tau)^{\gamma+3 a}\left\|\phi_x\right\|_{L^2}^2 d \tau \notag\\
= & \varepsilon \int_0^t(1+\tau)^\gamma\left\|\phi_{x x}\right\|_{L^2}^2 d \tau+C(1+t)^{\gamma+\frac{3}{2} b+2\delta+1}+\int_0^t(1+\tau)^{\gamma+3 a-3} {}(1+\tau)^3\left\|\phi_x\right\|_{L^2}^2 d \tau \notag\\
\leq & \varepsilon \int_0^t(1+\tau)^\gamma\left\|\phi_{x x}\right\|_{L^2}^2 d \tau+C(1+t)^{\gamma+\frac{3}{2} b+2\delta+1}+C(1+t)^{\gamma+3 a-3} {} \int_0^t(1+\tau)^3\left\|\phi_x\right\|_{L^2}^2 d \tau \notag\\
\leq & \varepsilon \int_0^t(1+\tau)^\gamma\left\|\phi_{x x}\right\|_{L^2}^2 d \tau+C(1+t)^{\gamma+\frac{3}{2} b+2\delta+1}+C(1+t)^{\gamma+3 a-3} {}(1+t)^{3-\frac{1}{2}+\frac{5}{4} \delta} \notag\\
\leq & \varepsilon \int_0^t(1+\tau)^\gamma\left\|\phi_{x x}\right\|_{L^2}^2 d \tau+C(1+t)^{\gamma+\left(\frac{3}{2} b+1\right)+2 \delta}+C(1+t)^{\gamma+\left(3 a-\frac{1}{2}\right)+2 \delta},
\end{align}
\noindent where $a,b$ satisfy $a+b=-\frac{4}{3}$. Now take $3a-\frac{1}{2}=\frac{3}{2}b+1$, thus $a=-\frac{1}{9}, b=-\frac{11}{9}$, which leads to 
\begin{equation}
     \label{J14 esti}
     \int_0^t(1+\tau)^\gamma \int_{\mathbb{R}}\left|\phi_x\right| {}|F_x| d x d \tau \leq \varepsilon \int_0^t(1+\tau)^\gamma\left\|\phi_{x x}(\tau)\right\|_{L^2}^2 d \tau+C(1+t)^{\gamma-\frac{5}{6}+2\delta}.
\end{equation}

\noindent Since $\delta>0$ is arbitrary, we may take it sufficiently small and then rename it. Substituting  (\ref{J6 esti}), (\ref{J7 esti}), (\ref{J8,J11,J12 esti}), (\ref{J9 esti}), (\ref{J10,J15,J16 esti}), (\ref{J14 esti}) and (\ref{J13 esti}) into (\ref{derivative weighted esti initial: degenerate case}), it follows that
\begin{equation}
  \begin{aligned}
& (1+t)^\gamma\left(\left\|\phi_x(t)\right\|_{L^2}^2+\left\|\phi_{t x}(t)\right\|_{L^2}^2+\left\|\phi_{x x}(t)\right\|_{L^2}^2\right)+\int_0^t(1+\tau)^\gamma\left(\left\|\phi_{t x}(\tau)\right\|_{L^2}^2+\left\|\phi_{x x}(\tau)\right\|_{L^2}^2\right) d \tau \\
\leq & C\left(\phi_0, \phi_1\right)+C(1+t)^{\gamma-1+\frac{5}{4} \delta}+C(1+t)^{\gamma-\frac{3}{2}+\frac{5}{4} \delta}+C(1+t)^{\gamma-\frac{5}{6}+2 \delta}+C(1+t)^{\gamma-\frac{13}{14}+\delta} \\
\leq & C(1+t)^{\gamma-\frac{5}{6}+2 \delta}.
\end{aligned}
\end{equation}
\noindent Thus we complete the proof of Lemma \ref{derivative L2 decay: degenerate case}. The Sobolev inequality, Lemma \ref{L2 decay: degenerate case}, and Lemma \ref{derivative L2 decay: degenerate case} yield
\begin{equation*}
\|\phi(t)\|_{L^{\infty}}\leq C\|\phi(t)\|_{L^2}^{1/2}\|\phi_x(t)\|_{L^2}^{1/2}
\leq C(1+t)^{-\frac{1}{3}+\delta}.
\end{equation*}
For $\psi=v-f(\widetilde U)$, the second equation of the relaxation system gives
\begin{equation*}
\psi_t+\psi=f(\widetilde U+\phi)-f(\widetilde U)-\phi_x
-\widetilde U_x-f'(\widetilde U)\widetilde U_t.
\end{equation*}
By the derivative estimate just proved and the Sobolev inequality,
$\|\phi_x(t)\|_{L^\infty}\leq C(1+t)^{-5/12+\delta}$.
The profile estimates also give
$\|\widetilde U_x(t)\|_{L^\infty}+\|\widetilde U_t(t)\|_{L^\infty}\leq C(1+t)^{-1/2}$.
Using the exponential representation of this damping equation, we obtain
$\|\psi(t)\|_{L^\infty}\leq C(1+t)^{-1/3+\delta}$.

Finally, we give the proof of Theorem \ref{Main Th 2.1}. Since for any \(\delta>0\), \(\|\phi(t)\|_{L^{\infty}}\leq C(1+t)^{-\frac{1}{3}+\delta}\), it follows from the uniform bound of \(\|\phi(t)\|_{L^{\infty}}\) that the inequality holds for \(t\geq 0\). Notice that for \(t\geq 1\),
    \begin{align}
\left\|\widetilde{U}-U_{\text {multi }}(t)\right\|_{L^{\infty}} & =\left\|U^r(t, x)+U(1+t, x)-\left(u^r\left(\frac{x}{t}\right)+U(t, x)\right)\right\|_{L^{\infty}} \notag\\
& \leq\left\|U^r(t, x)-u^r\left(\frac{x}{t}\right)\right\|_{L^{\infty}}+\|U(1+t, x)-U(t, x)\|_{L^{\infty}} \notag\\
& \leq C(1+t)^{-1+\frac{1}{4}}+\left\|\int_t^{1+t} U_s(s, x) d s\right\|_{L_x^{\infty}} \notag\\
& \leq C(1+t)^{-1+\frac{1}{4}}+\int_t^{1+t}\left\|U_s(s)\right\|_{L^{\infty}} d s \notag\\
& \leq C(1+t)^{-1+\frac{1}{4}}+C t^{-\frac{1}{2}} \int_t^{1+t} d s \leq C(1+t)^{-\frac{1}{2}},
\end{align}
\noindent which yields that for $\delta>0$ and $t\geq1$,
\begin{equation}
    \left\|u-U_{\text {multi }}(t)\right\|_{L^{\infty}} \leq\|\phi(t)\|_{L^{\infty}}+\left\|\widetilde{U}-U_{\text {multi }}(t)\right\|_{L^\infty} \leq C(1+t)^{-\frac{1}{3}+\delta}.
\end{equation}
\noindent Regarding the $L^1$ bound, since for $t>0$,
\begin{equation}
    \begin{aligned}
\left\|\widetilde{U}-U_{\text {multi }}(t)\right\|_{L^1} & \leq\left\|U^r(t, x)-u^r\left(\frac{x}{t}\right)\right\|_{L^1}+\|U(1+t, x)-U(t, x)\|_{L^1} \\
& \leq C+\left\|\int_t^{1+t} U_s(s, x) d s\right\|_{L_x^1} \\
& \leq C+\int_t^{1+t}\|U_s(s)\|_{L^1}\,ds \leq C,
\end{aligned}
\end{equation}
\noindent where we have used (5) of Lemma \ref{Properties of the smooth rarefaction wave 2}. Therefore, one can deduce that for \(\delta>0\) and \(t>0\),
\begin{equation}
    \left\|u-U_{\mathrm{multi}}(t)\right\|_{L^1} \leq\|\phi(t)\|_{L^1}+\left\|\widetilde{U}-U_{\mathrm{multi}}(t)\right\|_{L^1} \leq C(1+t)^\delta.
\end{equation}
\noindent Moreover, the Lipschitz continuity of $f$ on the state interval transfers the perturbation estimates to $v-f(U_{\mathrm{multi}})$. This concludes the proof of Theorem \ref{Main Th 2.1}. 

\bigskip
\noindent\textbf{Competing interests.} The authors declare that they have no competing interests.

\medskip
\noindent\textbf{Data availability.} Data sharing is not applicable to this article as no datasets were generated or analyzed during the current study.

\IfFileExists{actams_two_wave.bst}{\bibliographystyle{actams_two_wave}}{\bibliographystyle{amsplain}}
\bibliography{reference_two_wave}

\end{document}